\documentclass[10pt, a4paper]{article}

\usepackage{graphicx}
\usepackage{color}
\usepackage{amsmath}
\usepackage{amssymb}
\usepackage{amscd}
\usepackage{bbm}
\usepackage{tikz}
\usetikzlibrary{patterns}
\usepackage{hyperref}
\hypersetup{
    bookmarks=true,         
    colorlinks=true,       
    linkcolor=red,          
    citecolor=green,        
    filecolor=magenta,      
    urlcolor=cyan           
}

\usepackage[utf8]{inputenc}
\usepackage{amsthm}
\usepackage{bbm} 
\usepackage[linesnumbered,lined,boxed,commentsnumbered]{algorithm2e}

\usepackage{marginnote}

\newtheorem{Theorem}{Theorem}
\newtheorem{Assumption}{Assumption}
\newtheorem{Definition}{Definition}

\newtheorem{Proposition}{Proposition}

\newtheorem{Remark}{Remark}

\newcommand{\inr}[1]{\bigl< #1 \bigr>}

\newcommand{\norm}[1]{\left\|#1\right\|}%

\newcommand\eps{\epsilon}

\def \endproof
{{\mbox{}\nolinebreak\hfill\rule{2mm}{2mm}\par\medbreak}}

\DeclareMathOperator*{\argmin}{argmin}

\def\ds1{\textrm{1\kern-0.25emI}} 

\newcommand \cF{{\cal F}}

\newcommand \cI{{\cal I}}

\newcommand \cK{{\cal K}}

\newcommand \cN{{\cal N}}
\newcommand \cO{{\cal O}}

\newcommand \cS{{\cal S}}

\newcommand \bE{{\mathbb E}}

\newcommand \bN{{\mathbb N}}

\newcommand \bP{{\mathbb P}}

\newcommand \bR{{\mathbb R}}

\DeclareMathOperator*{\Med}{Med}
\DeclareMathOperator*{\Tr}{Tr}

\usepackage{braket}
\usepackage{authblk}

\begin{document}

\title{On the robustness to adversarial corruption and to heavy-tailed data of the Stahel-Donoho median of means}
\author[1]{Jules Depersin and Guillaume Lecu{\'e}  \\ email: \href{mailto:jules.depersin@ensae.fr}{jules.depersin@ensae.fr}, email: \href{mailto:lecueguillaume@gmail.com}{guillaume.lecue@ensae.fr} \\ CREST, ENSAE, IPParis. 5, avenue Henry Le Chatelier, 91120 Palaiseau, France.}

\date{}                     
\setcounter{Maxaffil}{0}
\renewcommand\Affilfont{\itshape\small}

\maketitle

\begin{abstract}We consider median of means (MOM) versions  of  the Stahel-Donoho outlyingness (SDO) \cite{stahel1981robuste,donoho1982} and of Median Absolute Deviation (MAD) \cite{MR0359096} functions to construct subgaussian estimators of a mean vector under adversarial contamination and heavy-tailed data. We develop a single analysis of the MOM version of the SDO which covers all cases ranging from the Gaussian case to the $L_2$ case. It is based on isomorphic and almost isometric properties of the MOM versions of SDO and MAD. This analysis also covers cases where the mean does not even exist but a location parameter does; in those cases we still recover the same subgaussian rates and the same price for adversarial contamination even though there is not even a first moment. These properties are achieved by the classical SDO median and are therefore the first non-asymptotic statistical bounds on the Stahel-Donoho median complementing the $\sqrt{n}$-consistency \cite{maronna1995behavior} and asymptotic normality \cite{MR2051003} of  the Stahel-Donoho estimators.  We also show that the MOM version of MAD can be used to construct an estimator of the covariance matrix under only a $L_2$-moment assumption or of a scale parameter if a second moment does not exist.
\end{abstract}


\section{Introduction } 
\label{sec:introduction}
Robust estimation of a mean vector has witnessed an important renewal during the last decade. Two communities have looked at this problem from their own perspective. In the statistics literature, several works have considered the problem of robustness with respect to heavy-tailed data. The aim here is to construct an estimator achieving statistical bounds with the same confidence as if all the data were i.i.d. Gaussian even  though the data at hand are only assumed to have a second moment. Such estimators are called subgaussian estimators; they are said to be robust to heavy-tailed data. The first seminal result showing the existence of such an estimator may be found in \cite{MR3052407}. It is also shown in \cite{MR3052407} that the empirical mean does not achieve this goal: the rate achieved by the empirical mean in the $L_2$ setup cannot be better than $\sigma\sqrt{1/(\delta N)}$ with probability at least $1-\delta$ whereas it is of the order of $\sigma\sqrt{\log(1/\delta)/N}$ when the data are i.i.d. Gaussian. The rate $\sigma\sqrt{\log(1/\delta)/N}$ is called the subgaussian rate for the mean estimation problem in the one-dimensional case. It is this rate that has been first achieved in \cite{MR3052407} by an $M$-estimator with a specific  score function. This rate was then achieved using a median-of-means principle in several works such as \cite{devroye2016,MR3124669}. It was then extended to the $d$-dimensional case in many other works \cite{lugosi2019sub,LMSL,MR3845006,hopkins2018sub,Bartlett19,Jules_Guillaume_1,lei2020fast,dalalyan2020all} since then.

For the mean estimation problem in $\bR^d$, most of the results have been given w.r.t. the Euclidean $\ell_2^d$ distance. There is however no statistical justification for this choice but that the $\ell_2^d$ metric is simply the most natural Hilbert metric in $\bR^d$ and so it seems natural to use it as a way to measure the statistical performance of an estimator of a $d$-dimensional vector. The resulting confidence sets have therefore the form $\hat \mu + r^*_{N,\delta}B_2^d$ where $\hat \mu$ is an estimator, $B_2^d=\{x\in\bR^d:\norm{x}_2\leq 1\}$ is the unit Euclidean ball and $r^*_{N,\delta}$ is the rate of convergence w.r.t. $\ell_2^d$ achieved by $\hat \mu$ with confidence $1-\delta$. When estimating w.r.t. the $\ell_2^d$ metric, confidence sets are therefore $\ell_2^d$-balls. One may wonder if these confidence sets are the best from a statistical point of view, for instance, the one with smallest volume for a fixed confidence $1-\delta$. To answer this type of question, one usually go back to the ideal i.i.d. Gaussian case, and use results obtained in that framework as benchmark results. We may also consider this model to design optimal benchmark confidence sets, that could be used to define more appealing estimation metric of a mean vector in $\bR^d$. 

Let us now see what are the ''best'' (in some sense given later) confidence sets in the i.i.d. Gaussian case: let $X_1,\ldots, X_N$ be i.i.d. distributed like $\cN(\mu, \Sigma)$ where $\mu\in\bR^d$ is the mean and $\Sigma$ is a symmetric definite positive  matrix (we assume here that $\Sigma$ is invertible). The MLE is the empirical mean $\bar X_N$ and $\sqrt{N}(\bar X_N-\mu)\sim \cN(0,\Sigma)$. The later result holds asymptotically if the data are only assumed to be in $L_2$ thanks to the CLT. The key observation here is that $\Sigma$ is the inverse of the Fisher information in this model and thus there are no regular asymptotically normal $M$-estimator that can estimate the mean with an asymptotic covariance matrix better than $\Sigma$. Moreover, level sets of the standard Gaussian density function are Euclidean $B_2^d$ balls centered at zero. As a consequence, the best confidence sets for $\mu$ with confidence $1-\delta$ are ellipsoids $\Sigma^{-1/2}B_2^d$ with radius given by the quantile of order $1-\delta$ of a chi-square variable with parameter $d$ centered at the estimator. This type of confidence region are equivalently written as estimation results of $\mu$ with respect to the norm $x\in\bR^d\to \norm{\Sigma^{-1/2}x}_2$. It follows that the best metric -- that is the one leading to minimal volume confidence sets for a given confidence in the benchmark i.i.d. Gaussian case  -- is the norm $\norm{\Sigma^{-1/2}\cdot}_2$ whose unit ball is the ellipsoid  $\Sigma^{1/2}B_2^d$. 

Regarding our robust mean estimation problem, the two next natural questions are the following: is it possible to construct robust mean estimators w.r.t. the $\norm{\Sigma^{-1/2}\cdot}_2$ metric? and what is the best convergence rate one can hope for? In the literature \cite{MR4026610,convex_prog}, one may find estimators which can estimate in a robust way a mean vector w.r.t. any metric of the type $u\in\bR^d\to\norm{u}_S=\sup_{v\in S}\inr{v,u}$ where $S\subset\bR^d$. In particular, for $S=\Sigma^{1/2}B_2^d$, this metric coincides with the one we want to use, i.e. $\norm{\Sigma^{-1/2}\cdot}_2$. It has also been proved that the optimal deviation minimax rate (the one obtained in the benchmark i.i.d. Gaussian case) is for the mean estimation problem with respect to $\norm{\cdot}_S$ given by (see \cite{convex_prog}) 
\begin{equation}\label{eq:subgaussian_rate}
\sqrt{\frac{\ell^*(\Sigma^{1/2}S)}{N}} + \sup_{v\in S}\norm{\Sigma^{1/2}v}_{2}\sqrt{\frac{\log(1/\delta)}{N}}.
\end{equation}For instance, for $S=B_2^d$ that is for $\norm{\cdot}_S = \norm{\cdot}_2$, the later rate is the classical $\sqrt{\Tr(\Sigma)/N} + \sqrt{\norm{\Sigma}_{op}\log(1/\delta)/N}$ rate. The case that is interesting to us is when $\norm{\cdot}_S = \norm{\Sigma^{-1/2}\cdot}_2$, that is for $S = \Sigma^{-1/2}B_2^d$. In that case, the subgaussian rate is 
\begin{equation}\label{eq:subgaussian_rate_Sigma}
\sqrt{\frac{d}{N}} + \sqrt{\frac{\log(1/\delta)}{N}}.
\end{equation}This is the rate we will try to reach from an adversarial corrupted and heavy-tailed dataset. We will also have to take into account the price for corruption. There are indeed known information theoretic lower bounds showing that there are no statistics that can do better than $(|\cO|/N)^\alpha$ where $\alpha\in[1/2,1]$ is some exponent depending on properties of the good data. For instance, $\alpha=1$ for Gaussian variables and $\alpha=1/2$ for some $L_2$ variables. However, we will see that the best possible cost $|\cO|/N$ (i.e. for $\alpha=1$) can be achieved even variables which do not have a first moment as long as the cdfs of all one-dimensional projections of the centered and normalized data are regular enough. 

Unfortunately, all estimators known to achieve the subgaussian rate in \eqref{eq:subgaussian_rate_Sigma} (the Le Cam test estimator in \cite{MR4026610}, the minmax MOM estimator with loss function $\ell(x,u)=\norm{\Sigma^{-1/2}(u-v)}$ from \cite{lerasle2019monk} or the Fenchel-Legendre estimators from \cite{convex_prog}) are using the set $S$ in their construction. This is something we cannot do here because $S=\Sigma^{-1/2}B_2^d$ depends on $\Sigma$ which is unknown in general. One therefore has to consider other type of estimators than the ones cited above. In this work,  we will do it thanks to a notion of depth/outlyingness introduced at the beginning of the $80$'s which, unlike the last cited estimators, uses a normalization by a robust estimation of the scale.

There are several ways to measure how 'deep' is a vector with respect to a cloud of points, see for instance the half-space depth of Tukey \cite{MR0426989,nagy2019halfspace}, the simplicial depth \cite{MR1041400,MR1186260}, Mahalanobis depth or the projection depth \cite{MR1214839}. Taking a point with maximal depth is usually seen as a way to define a median in $\bR^d$ (see Radon points \cite{barany2020application} or Fermat Points \cite{haldane1948note}). There are therefore several ways to define a median of a cloud of points in $\bR^d$. One depth has received a particular attention both in theory and in practice and is known as the Stahel-Donoho  outlyingness (SDO) \cite{stahel1981robuste,donoho1982breakdown}. It can be used to construct estimators of multivariate location and scatter known as the Stahel-Donoho estimators (SDE) which were the first equivariant estimators with a high breakdown point. The aim of this work is to show that this notion of depth can be used to construct estimator of a mean vector in $\bR^d$ which is robust to adversarial contamination and to heavy-tailed data with respect to $\norm{\Sigma^{-1/2}\cdot}_2$. Let us now define this notion of depth\footnote{We speak here about depth instead of outlyningness : these two concepts are expressing the same notion but in reverse order.} and recall some of its properties.

There is a common approach to many notion of depths for a general $d$-dimensional set of vectors: first, a definition of depth in $\bR$ is given and second, this notion is extended to $\bR^d$ simply by applying this one-dimensional definition to the set of one-dimensional projections of the data in all directions $v\in\bR^d$ (or all $v\in S$ for some subset $S\subset \bR^d$) and then by taking the supremum over all $v\in\bR^d$ (or $v\in S$). This approach is based on the idea that if a point in $\bR^d$ is an outlier then there must be some direction $v$ such that it is an (univariate) outlier when projected into that direction.

The SDO of $z\in\bR$ with respect to a dataset $\{a_1, \ldots, a_K\}$ in $\bR$ is defined as  
\begin{equation}\label{eq:def_SDO_intro}
SDO(z;\{a_1, \ldots, a_K\}) = \frac{|z-\Med(a_k)|}{\Med(|a_k - \Med(a_k)|)}
\end{equation}and a natural extension to $\bR^d$ is using the previous one for all one-dimensional projections of the data and by taking the supremum over all directions: for any $\nu\in\bR^d$ and a dataset $\{Z_1, \ldots, Z_K\}$ in $\bR^d$, we set
\begin{align}\label{eq:intro_SDO_Rd}
\notag SDO(\nu, \{Z_1, \ldots, Z_K\})& = \sup_{v\in\bR^d} SDO(\inr{\nu, v}; \{\inr{Z_1, v}, \ldots, \inr{Z_K, v}\})\\
&=\sup_{v\in\bR^d}\frac{|\inr{\nu, v}-\Med(\inr{Z_k, v})|}{\Med(|\inr{Z_k, v} - \Med(\inr{Z_k,v})|)}.
\end{align}A natural way to define a median of the $Z_k$'s is obtained by taking a point with minimal outlyingness (i.e. maximal depth): 
\begin{equation*}
\hat \mu^{SDO} \in\argmin_{\mu\in\bR^d} SDO(\mu, \{Z_1, \ldots, Z_K\}).
\end{equation*}However $ \hat \mu^{SDO} $ is not the most usual choice to estimate some location of the $Z_k$'s when they are assumed to follow a statistical model. The Stahel-Donoho location estimator is rather defined as a convex sums of the data:
\begin{equation}\label{eq:SD_estimator_location_intro}
\hat \mu^{SDE}_K = \frac{\sum_{k=1}^K w_k Z_k}{\sum_{k=1}^K w_k} 
\end{equation}where the weights are some function of the outlyingness of the data, i.e. $w_k = w(SDO(Z_k))$ for some (decreasing) weight function $w:\bR^+\to \bR^+$. The weights can also be used to estimate the scatter of the points by 
\begin{equation}\label{eq:SD_estimator_scatter_intro}
\hat \Sigma^{SDE}_K = \frac{\sum_{k} w_k (Z_k-\hat \mu^{SDE})(Z_k-\hat \mu^{SDE})^\top}{\sum_k w_k}. 
\end{equation} Note that there is a more general definition of SDO than the one considered in \eqref{eq:def_SDO_intro} with general (one dimensional) definitions of location and scale statistics; in \eqref{eq:def_SDO_intro}, we used the median $\Med(a_k)$ and Median Absolute Deviation (MAD) $\Med(|a_k - \Med(a_k)|)$ for these statistics \cite{MR362657}. 

As mentioned previously several results on the Stahel-Donoho Estimator (SDE) have been established during the last forty years. They are \textit{affine equivariant} meaning that for any affine transformation $x\in\bR^d\to A x +b$ of the dataset by a nonsingular matrix $A\in\bR^{d\times d}$ and a vector $b\in\bR^d$ the location estimator $\hat \mu^{SDE}_K$ is following the same transformation and the scatter estimator  $\hat \Sigma^{SDE}_K$ is transformed via $M\in\bR^{d\times d}\to A M A^\top$. SDE have been proved to have a \textit{finite-sample breakdown point} \cite{MR689745} which is the ''smallest amount of contamination necessary to upset an estimator entirely'' from \cite{MR1193313} in \cite{donoho1982}. In \cite{MR1292555}, it is proved that the SDE with MAD replaced by the average of the $k_1$th and $k_2$th smallest absolute deviation about the median $\Med(a_k)$ for $k_1=d-1+[(K+1)/2]$ and $k_2=d-1+[(K+2)/2]$ achieves the best finite-sample replacement breakdown point among all affine equivariant estimators obtained in \cite{MR902258} which is $[(K-d+1)/2]/K$ (this result holds when the weight function $w$ is continuous and there is an absolute constant $c_0$ such that $w(r)\leq c_0$, $w(r)\leq c_0/r^2$ for all $r\geq0$). This result was later extended in Theorem~3.2 from \cite{MR2051003}. This is to our knowledge the only established non-asymptotic properties of Stahel-Donoho estimators. 

There are however several asymptotic results for SDE such as a $\sqrt{n}$-consistency in \cite{maronna1995behavior}: if the $Z_k$'s are i.i.d. then $\sqrt{K}\left( (\hat \mu^{SDE}_K, \hat \Sigma^{SDE}_K) - (\mathbf{t}, \mathbf{V})\right)$ tends to $0$ in probability when $K\to+\infty$ where $\mathbf{t}$ and $\mathbf{V}$ are some location and scatter parameters of the distribution of $Z_1$. This result holds when the weight function $w$ is such as $|w(r)-w(r^\prime)|\leq \gamma \min(1,1/\min(r,r^\prime)^3)|r-r^\prime|$ for all $r,r^\prime\in\bR$ and when for all $v\in\bR^d$ the cumulative distribution function (cdf) of $\inr{Z_1,v}$ denoted by $F_v$ satisfies the following assumption : there exists some absolute constants $c_0>0$ and $c_1>0$ such that for all $|\eps|\leq c_0$
\begin{equation}\label{eq:ass_maronna_intro}
|F_v(\Med(F_v)+\eps) - F_v(\Med(F_v))|\geq c_1 |\eps| \mbox{ and } |F_v(\Med(F_v)\pm \sigma_v+\eps) - F_v(\Med(F_v\pm \sigma_v))|\geq c_1 |\eps|
\end{equation}where $\Med(F_v)=\inf(x\in\bR: F_v(x)\geq1/2)$ is the median of $F_v$ and $\sigma_v=\Med(G_v) $ where $G_v$ is the cumulative distribution of the random variable $MAD(\inr{Z_1,v}):=\Med(|\inr{Z_1,v} - \Med(\inr{Z_1,v})|)$. A typical situation mentioned in \cite{maronna1995behavior} where \eqref{eq:ass_maronna_intro} holds is when the cdf $F:\bR^d\to [0,1]$ of $Z$ is such that $F=(1-\eta)F_0 + \eta F^*$ where $\eta<1$ and $F^*$ is any cdf and $F_0$ is such that there exists $c_0>0$ and $c_1>0$ such that for all $v\in\bR^d$, $\inr{Z_1,v}$ has a density denoted by $f_v$ satisfying $f_v(t)\geq c_1$ for all $t\in[\Med(F_v)\pm c_0]\cup[\Med(F_v)-\sigma_v \pm c_0]\cup [\Med(F_v)+\sigma_v \pm c_0]$. According to \cite{maronna1995behavior}, the later holds when $F$ is spherical with positive density in a neighborhood of $0$ and $\sigma_{e_1} e_1$ where $e_1=(1, 0, \cdots, 0)\in\bR^d$. We will come back later on these conditions since we will encounter similar assumptions for our analysis. Finally, asymptotic normality of SDE location estimators have been obtained in \cite{MR2051003} under great generality for the location and scatter estimators as well as for the weight function including the median and MAD estimators as in \eqref{eq:def_SDO_intro} and the projection depth obtained for the weight function $w:r\in\bR^+\to 1/(1+r)$. From a stochastic point of view, asymptotic results for $\hat \mu^{SDE}_K$ hold when the cdf $F$ is \textit{elliptically symmetric around $\mu$} which means that there exists a symmetric definite positive matrix $\Sigma$ such that for all $v\in\cS_2^{d-1}:=\{v\in\bR^d:\norm{v}_2=1\}$, $\inr{\Sigma^{-1/2}(Z_1-\mu),v}$ has the same distribution as $\inr{\Sigma^{-1/2}(Z_1-\mu),e_1}$ which is a univariate symmetric variable with density function $f$. In that case, asymptotic normality was obtained when $f(0)f(\sigma)>0$ where $\sigma = MAD(\inr{\Sigma^{-1/2}(Z_1-\mu),e_1})$. Again we will meet this type of condition in our analysis.

On the practical side, SDEs have been used a lot in practice and implementation on various languages such as R exists; and that is one reason why the study of the SDO may be useful, maybe more than some other notions of depth.  In the original paper \cite{stahel1981robuste}, the author proposes a random algorithm where the supremum over all directions  $v\in\bR^d$ is approximated by subsampling orthogonal directions to $d-1$ hyperplanes generated by $d$ randomly chosen points in the dataset.  Other strategies mixing random and deterministic directions have been proposed for instance in \cite{pena2007combining}. Several adaptations and extensions of this algorithm may be found in \cite{MR2752147} for an extension to an arbitrary kernel space or in \cite{van2011stahel,van2016stahel} for a ''cell-wise weights'' extension of the SDO where each coordinate of each data receives its own weight. However, only very little is known on the theoretical computational side. In Section~5 of \cite{donoho1992breakdown}, an algorithm running in time $\cO(K^{d+1}\log K)$ is mentioned but its time complexity is making this approach impractical for dimensions larger than $5$. There are to our knowledge no theoretical result of any kind on the convergence of some approximate algorithm for the computation of the SDO of a point in $\bR^d$ that could be used in practice. As mentioned already  in \cite{donoho1992breakdown}, ''some sort of computational breakthrough is necessary to make the estimators, as defined here, really practical''. This looks to be still the case. We will however not discuss about this issue in the present work and leave this question still opened.

The aim of this work is to construct mean vector estimators robust to adversarial outliers and heavy-tailed data achieving the deviation-minimax subgaussian rate from \eqref{eq:subgaussian_rate_Sigma} with respect to the metric $\norm{\Sigma^{-1/2}\cdot}_2$. 
On our way to our goal, we complement the results on the $\sqrt{n}$-consistency and the asymptotic normality of SDE,  by deriving the first non-asymptotic convergence rate for the original SDO median (as well as its median of means version) as a robust mean estimator in $\bR^d$ under the following assumption.
\begin{Assumption}\label{assum:first}[Adversarial contamination and $L_2$ inliers] There exists $N$ random vectors $(\tilde X_i)_{i=1}^N$ in $\bR^d$ which are independent with mean $\mu$ and covariance matrix  $\Sigma$. The $N$ random vectors $(\tilde X_i)_{i=1}^N$ are first given to an ''adversary'' who is allowed to modify up to $|\cO|$ of these vectors. This modification does not have to follow any rule. Then, the ''adversary'' gives back the modified dataset $(X_i)_{i=1}^N$ to the statistician. Hence, the statistician receives an ''adversarially'' contaminated dataset of $N$ vectors in $\bR^d$ which can be partitioned into two groups: the modified data $(X_i)_{i\in\cO}$, which can be seen as outliers and the ''good data'' or inliers $(X_i)_{i\in\cI}$ such that $\forall i\in\cI, X_i=\tilde X_i$. Of course, the statistician does not know which data has been modified or not so that the partition $\cO\cup\cI=\{1, \ldots, N\}$ is unknown to the statistician.  
\end{Assumption}The contamination model defined in Assumption~\ref{assum:first} covers the Huber $\eps$-contamination model from \cite{MR2488795} and also the $\cO\cup\cI$ contamination framework from \cite{MR4102681}. It has been popularized by the Computer Science community in particular in \cite{MR3631028,diakonikolas2017being}.  In the adversarial contamination model from Assumption~\ref{assum:first}, the set  $\cO$ can depend arbitrarily on the initial data $(\tilde X_i)_{i=1}^N$;  the corrupted data $(X_i)_{i \in \cO}$ can have any arbitrary dependence structure; and the informative data $(X_i)_{i\in\cI}$ may also be correlated (for instance, it is the case, in general,  when the $|\cO|$ data $\tilde X_i$ with largest $\ell_2^d$-norm are modified by the adversary).

In the setup defined by Assumption~\ref{assum:first}, we will use the SDO as one of our building block to achieve our goal as well as the Median-of-means principle \cite{MR702836,MR1688610,MR855970}. This principle has been extensively used during the last decades in particular for the problem of robust mean estimation \cite{LO,MR3576558,minsker2015geometric,MR4026610,MR4017683,LMSL,convex_prog,hopkins2018sub,Bartlett19}. The starting point of MOM estimator is to chose an integer $K\in[N]$, split the dataset into $K$ equal size blocks $B_1\sqcup\cdots\sqcup B_K=[N]$ (w.l.o.g. we assume that $N$ can be divided by $K$) and construct $K$ empirical means $\bar X_k = |B_k|^{-1}\sum_{i\in B_k}X_i$, one over each block. The Stahel-Donoho Median-of-Means that will be used to achieve the subgaussian rate \eqref{eq:subgaussian_rate_Sigma} with respect to $\norm{\Sigma^{-1/2}\cdot}_2$ under in the adversarial and heavy-tailed setup from Assumption~\ref{assum:first} is 
\begin{equation*}
\hat \mu^{SDO}_{MOM,K} \in\argmin_{\mu\in\bR^d} \sup_{\norm{v}_2=1}\frac{|\inr{\mu, v}-\Med(\inr{\bar X_k, v})|}{\Med(|\inr{\bar X_k, v} - \Med(\inr{\bar X_k,v})|)}.
\end{equation*}

Unlike recently introduced robust mean estimators, $\hat \mu^{SDO}_{MOM,K}$ is using a robust scatter estimator for normalization. Here its is a MOM version of MAD which is used to construct $\hat \mu^{SDO}_{MOM,K}$, i.e. $v\to \Med(|\inr{\bar X_k, v} - \Med(\inr{\bar X_k,v})|)$.  We will show that this normalization plays a central role in the analysis when one wants results w.r.t. the $\norm{\Sigma^{-1/2}\cdot}_2$-norm. But beyond this observation, we will show that MAD and its MOM version  satisfy isomorphic and almost-isometric properties that can be used for other task such as to construct estimator of the covariance under only a $L_2$ assumption as in Section~\ref{sec:estimation_of_} below.

The paper is organized as follows. In the next section, we consider the case where the good data have a Gaussian distribution and the dataset has been adversarially corrupted. In that case, no need to construct bucketed means and the original Stahel-Donoho median is proved to achieve the subgaussian rate \eqref{eq:subgaussian_rate_Sigma}. The Section~\ref{sec:the_L2_case} considers the general adversarial corrupted and heavy-tailed framework from Assumption~\ref{assum:first} where the MOM version of the SDO is proved to achieve the subgaussian rate. We also exhibit in this section a family of cdfs denoted here by $(H_{N,K,v}:v\in\cS_2^{d-1})$ which plays a key role in our analysis. In particular, when the behavior of these function around $0$ is similar to the one described above in \eqref{eq:ass_maronna_intro} then the same result as in the Gaussian case can be obtained and that may hold without anymore than $2$ moments (see Section~\ref{sec:the_l2_case_via_behavior_in_zero}). In Section~\ref{sec:estimation_of_}, we show how to use the MOM version of MAD to construct an estimator of the covariance matrix under $L_2$. In Section~\ref{sec:study_of_the_}, we explore the properties of the family of functions  $(H_{N,K,v}:v\in\cS_2^{d-1})$. A conclusion and open questions are provided in Section~\ref{sec:conclusion} that are followed by the proofs of all the results in Section~\ref{sec:proofs}.

\paragraph{Notations.} We denote by $x\in\bR^d\to\norm{x}_2=\left(\sum_j x_j^2\right)^{1/2}$ the Euclidean norm with associated unit sphere $\cS_2^{d-1}$ and ball $B_2^d$. We also denote by $g\sim\cN(0,1)$ a standard one-dimensional Gaussian variable and its associated standard Gaussian cdf by $\Phi:t\in\bR\to \bP[g\leq t]=\int_{-\infty}^t \phi(u)du$ where $\phi:u\in\bR\to (2\pi)^{-1/2}\exp(-u^2/2)$ is the one dimensional Gaussian density function. We also set $H_G:t\to 1-\Phi(t)$ and $W_G:p\in(0,1)\to H_G^{(-1)}(p)$ the inverse function of $H_G$ so that $W(p)=\Phi^{-1}(1-p)$.

\section{The Gaussian case} 
\label{sec:the_gaussian_case}
In this section, we prove that the original SDO median achieves the (non-asymptotic) subgaussian rate \eqref{eq:subgaussian_rate_Sigma} when the dataset may have been corrupted by an adversary and when the good data have a Gaussian distribution; our main model assumption is the following.
\begin{Assumption}\label{assum:gaussian}[Adversarial contamination and Gaussian inliers] There exists $N$ i.i.d. Gaussian vectors $(G_i)_{i=1}^N$ in $\bR^d$ with mean $\mu$ and (unknown) covariance matrix $\Sigma$. We assume that $\Sigma$ is invertible. The $N$ random vectors $(G_i)_{i=1}^N$ are first given to an ''adversary'' who is allowed to modify up to $|\cO|$ of these vectors. This modification does not have to follow any rule. Then, the ''adversary'' gives the modified dataset $(X_i)_{i=1}^N$ to the statistician. 
\end{Assumption}
We use the Gaussian case again as a benchmark case for the more involved heavy-tailed situation which requires to bucket the data and some assumption on the distribution of the good data. When the ''good'' data are Gaussian there  is no need to bucket the data and the elliptically symmetric property of the Gaussian variables is simplifying the analysis. The mean estimator we use in this section is therefore the median of the original  Stahel-Donoho outlyingness function
\begin{equation}\label{eq:SDO_Gauss}
SDO: \mu\in\bR^d \to  \sup_{v\in\bR^d}\frac{|\inr{\mu, v}-\Med(\inr{X_i, v})|}{\Med(|\inr{X_i, v} - \Med(\inr{X_i,v})|)}
\end{equation}and the associated median is a point minimizing this outlyingness function:
\begin{equation*}
\hat\mu^{SDO}\in \argmin_{\mu\in\bR^d} SDO(\mu).
\end{equation*}

Our main result in the adversarial corruption setup with Gaussian inliers is the following:

\begin{Theorem}\label{thm:SDO_gauss_robust_subgauss} There are absolute constants $c_0,c_1$ and $c_2$ such that the following holds.
We assume that the adversarial contamination with Gaussian inliers model Assumption~\ref{assum:gaussian} holds with a number of adversarial outliers denoted by $|\cO|$. Let $0<\eps<\Phi^{-1}(3/4)/c_1$.  We assume that $|\cO|\leq \eps N$ and $N\geq c_0 \eps^{-2}d$. For all $0<u<c_1\eps^2 N$, with probability at least $1-2\exp(-u)$, 
\begin{equation*}
\norm{\Sigma^{-1/2}(\hat \mu^{SDO} - \mu)}_2\leq 2\left(1 + c_1 \eps\right)\left(C_0\sqrt{\frac{d+1}{N}} + \sqrt{\frac{u}{N}} + \frac{|\cO|}{N}\right).
\end{equation*}
\end{Theorem}
Let us first remark that if $N<d$ then the $N$ data $X_1, \ldots, X_N$ cannot span the entire $\bR^d$ space and so there exists a non zero vector $v\in\bR^d$ which is orthogonal to all the data points. Hence, $MAD(v):=\Med(|\inr{X_i, v} - \Med(\inr{X_i,v})|) = 0$ a.s. and so $SDO(\mu)=+\infty$ for all $\mu\in\bR^d$. Therefore, assuming that $N\geq d$ is a minimal assumption when we work with the SDO function. We also note that the factor $\Phi^{-1}(3/4)$ is sometimes used as a renormalization factor in the definition of the Stahel-Donoho outlyingness function \cite{rousseeuw2018measure}.

Theorem~\ref{thm:SDO_gauss_robust_subgauss} shows that the SD median $\hat \mu^{SDO}$ is robust to adversarial contamination up to a proportion of $N$ and that the rate achieved remains the same as if there was no contamination when $|\cO|\lesssim \sqrt{N}\min(\sqrt{u}, \sqrt{d})$. If we put this result with regard to the finite-sample replacement breakdown point (RBP) achieved by the SDE with a slight modification of MAD at the denominator as recalled in the Introduction, we see that the order of magnitude are the same: SDE and $\hat \mu^{SDO}$ can handle both a proportion of $N$ adversarial outliers. The constant in the RBP (which is close to $1/2$ when $N>>d$) is certainly better than the one obtained in T{}heorem~\ref{thm:SDO_gauss_robust_subgauss} but the result in the later theorem shows that the estimator still achieve the deviation minimax rate \eqref{eq:subgaussian_rate_Sigma} even up to $\eps N$ outliers whereas RBP can only insure that the estimator does not go to infinity; RBP does not ensure any statistical convergence rate after data corruption unlike Theorem~\ref{thm:SDO_gauss_robust_subgauss} does. 

The rate of convergence obtained in Theorem~\ref{thm:SDO_gauss_robust_subgauss} has been obtained by several other procedures. For instance, it has been proved that the Tukey median achieves this rate in \cite{MR3845006} when the covariance is proportional to the identity and for the Huber-contamination setup. The same bound was also obtained by a polynomial time algorithm in \cite{dalalyan2019outlier} when the covariance matrix $\Sigma$ is known.

The proof of Theorem~\ref{thm:SDO_gauss_robust_subgauss} (which may be found in Section~\ref{sec:proofs}) is based on two isomorphic principles of the MAD and SDO functions. We will extend these two properties to the MOM versions of MAD and SDO in the next section. For the moment, let us recall their definitions and write these two properties that are interesting beyond the proof of Theorem~\ref{thm:SDO_gauss_robust_subgauss}.

The normalization factor in the SDO function \eqref{eq:SDO_Gauss} is called the MAD (median absolute deviation) \cite{MR0359096}
\begin{equation*}
MAD:v\in\bR^d\to \Med(|\inr{X_i, v} - \Med(\inr{X_i,v})|).
\end{equation*}It plays a key role to get estimation result w.r.t. the $\norm{\Sigma^{-1/2}\cdot}_2$ norm whereas $\Sigma$ is unknown. However, this normalization factor requires some more work than for the analysis of classical robust estimators that are only focused on the estimation of the mean. Indeed, $MAD(v)$ is actually a robust estimator of the scatter of $\inr{G,v}$ which is $ \Phi^{-1}(3/4)\norm{\Sigma^{1/2}v}_2$ (note that if $g\sim\cN(0,1)$ then $MAD(g)=\Phi^{-1}(3/4)$). It is therefore a 'second order' robust estimator but since it appears in the denominator of the SDO function, we cannot only prove an upper estimate for this quantity and we need an isomorphic result -- that is upper and lower matching (up to constants) bounds -- on the MAD. This result is of independent interest and we are therefore stating it here. The proof is given in Section~\ref{sub:proof_of_the_props}. We also state a similar isomophic result for SDO which can be use to prove Theorem~\ref{thm:SDO_gauss_robust_subgauss}. We will see later in Section~\ref{sec:study_of_the_} that these metric properties of SDO and MAD can be extended to cases where the mean does not even exists (in that case $\mu$ is a \textit{location} parameter) showing that these properties have actually more to do with elliptical symmetry than they have to do with concentration.

\begin{Proposition}\label{prop:MAD} There are absolute constants $c_1,c_2$ and $c_3$ such that the following holds. Let $0<\eps<\Phi^{-1}(3/4)/c_1$.
We assume that the adversarial model with Gaussian inliers Assumption~\ref{assum:gaussian} holds with a number of adversarial outliers $|\cO|\leq \eps N$. We assume that $N\geq c_2\eps^{-2} d$. With probability at least $1-\exp(-c_3 \eps^2N)$, for all $v\in\bR^d$,
\begin{equation*}
(\Phi^{-1}(3/4) - c_1 \eps)\norm{\Sigma^{1/2}v}_2\leq MAD(v)  \leq (\Phi^{-1}(3/4) + c_1 \eps)\norm{\Sigma^{1/2}v}_2.
\end{equation*}

Moreover, for all $0<u<c_1\eps^2 N$, with probability at least $1-2\exp(-u)$, for all $v\in\bR^d$,  if $\norm{\Sigma^{-1/2}(v-\mu)}_2\geq 2r^*$ then
\begin{equation*}
\frac{\norm{\Sigma^{-1/2}(v-\mu)}_2}{2(\Phi^{-1}(3/4) + c_1 \eps)\sqrt{K/N}}  \leq SDO_K(v)\leq \frac{3\norm{\Sigma^{-1/2}(v-\mu)}_2}{2(\Phi^{-1}(3/4) - c_1 \eps)\sqrt{K/N}} 
\end{equation*}and if $\norm{\Sigma^{-1/2}(v-\mu)}_2\leq 2r^*$ then $ SDO_K(v)\leq 3 r^*(\Phi^{-1}(3/4) - c_1 \eps)^{-1}$ where $\Phi^{-1}$ is the quantile function of a standard Gaussian cdf and $r^*=\left(C_0\sqrt{(d+1)/N} + \sqrt{u/N} + |\cO|/N\right)$ is the subgaussian rate from \eqref{eq:subgaussian_rate_Sigma} with the additive adversarial contamination term $|\cO|/N$.
\end{Proposition}

The isomorphic properties of the MAD and SDO functions uniformly over $\bR^d$ imply the robustness and subgaussian properties of the SDO median in Theorem~\ref{thm:SDO_gauss_robust_subgauss}. Similar results for other depths may be found in the literature on robust mean estimation such as the isomorphic property of the Tukey depth proved in \cite{MR3845006}.

\section{The $L_2$ case} 
\label{sec:the_L2_case}
In this section, we do not anymore assume that the good data follow a Gaussian distribution but we only assume that they have a second moment (and that the dataset may still be contaminated by an adversary following Assumption~\ref{assum:first} ). Nevertheless, even though we are in the heavy tailed setup with adversarially corrupted data we still want to achieve the subgaussian rate for the $\norm{\Sigma^{-1/2}\cdot}_2$-norm. To achieve such a result the median-of-means principle has been proved to perform well. We will therefore use this principle together with the Stahel-Donoho concept of outlyingness. We introduce now an estimator constructed according to these two principles.

Let $K\in[N]$ be the number of blocks and let $\bar X_k=(1/|B_k|)\sum_{i\in B_k} X_i, k\in[K]$ be the bucketed means.  Outlyingness / depth of a point $\mu\in\bR^d$ is measured with respect to the bucketed means:
\begin{equation*}
SDO_K(\mu) = \sup_{v\in\bR^d}\frac{|\inr{\mu, v}-\Med(\inr{\bar X_k, v})|}{\Med(|\inr{\bar X_k, v} - \Med(\inr{\bar X_k,v})|)}
\end{equation*}and the Stahel-Donoho Median of means is defined as
\begin{equation*}
\hat \mu^{SDO}_{MOM,K} \in\argmin_{\mu\in\bR^d} SDO_K(\mu).
\end{equation*}

As for the Gaussian case, the isomorphic and nearly-isometric properties of $SDO_K$ and its denominator, called $MOMAD_K$, play a key role in our analysis. The $MOMAD_K$ is a Median of means version of the Median Absolute Deviation function. We denote it as MOMAD for Median Of Means Absolute Deviation:
\begin{equation}\label{eq:MOMAD}
 MOMAD_K:v\in\bR^d\to \Med\left(|\inr{\bar X_k, v} - \Med(\inr{\bar X_k,v})|\right).
 \end{equation} In the next section, we study metric properties of $MOMAD_K$  and for $SDO_K$ that will useful for our analysis of $\hat \mu^{SDO}_{MOM,K}$. Then, we will turn to the statistical bounds obtained for the median $\hat \mu^{SDO}_{MOM,K}$ in the  general heavy-tailed $L_2$ setup in Section~\ref{sec:the_L2_case_via_markov} and then we will study some extra regularity assumption of the cdfs $(H_{N,K,v}:v\in\cS_2^{d-1})$ at $0$ that allows to get better rates in Section~\ref{sec:the_l2_case_via_behavior_in_zero}.

\subsection{Some isomorphic and almost isometric properties of $MOMAD_K$ and $SDO_K$} 
\label{sub:an_isomorphic_property_of_}
In this section, we show that the MOM versions of the SDO and MAD operators (called $SDO_K$ and $MOMAD_K$) satisfy an isomorphic and almost-isometry properties as the $MAD$ and $SDO$ do in Proposition~\ref{prop:MAD} that holds only under a $L_2$ moment assumption.

 We introduce two families of functions which play a central role in our analysis. They involve the non-corrupted random variables $\tilde X_i, i\in[N]$ (and not the corrupted data $X_i,i\in[N]$).
 \begin{Definition}\label{def:H_and_W}
  For all $v\in\cS_2^{d-1}$,
\begin{equation}\label{eq:cdf_quantile}
H_v:=H_{N,K,v}:r\in\bR\to \bP\left[\frac{1}{\sqrt{N/K}}\sum_{i=1}^{N/K} \inr{\Sigma^{-1/2}(\tilde X_i-\mu), v}\geq r\right] \mbox{ and } W_v:=W_{N,K,v}:p\in(0,1)\to H_v^{(-1)}(p),
\end{equation}where $H_v^{(-1)}(p)=\max(r\in\bR: H_v(r)\geq p)$ is the generalized inverse of $H_v$.
 \end{Definition}

 As already observed in the proof of the $\sqrt{n}$-consistency of SDE from \cite{maronna1995behavior} as well as its asymptotic normality in \cite{MR2051003}, the behavior of the one-dimensional projection cdfs at the median and the two $1/4$ and $3/4$ quartiles play a central role in the analysis of SDO based estimators. This will also be the case for the MOM version of the SD median. It will appear in Section~\ref{sec:study_of_the_} that taking bucketed mean may force toward the Gaussian case for which all these conditions are naturally satisfied because of the elliptical symmetry of Gaussian variables. Let us now state our main assumption on the behavior of the one-dimensional quantile functions $(W_v:v\in\cS_2^{d-1})$. 

 \begin{Assumption}\label{ass:quantile_fct_empirical_mean}
 There exists some $0<\eps<1/8$ and some absolute constants $0<\varphi_l(\eps)<\varphi_u(\eps)$ such that for all $v\in\cS_2^{d-1}$,
\begin{equation*}
\max\left(W_v\left(\frac{1}{4}-2\eps\right) - W_v\left(\frac{1}{2}+2\eps\right), W_v\left(\frac{1}{2}-2\eps\right) - W_v\left(\frac{3}{4}+2\eps\right)\right)\leq \varphi_u(\eps)
\end{equation*}and
\begin{equation*}
\min\left(W_v\left(\frac{1}{4}+2\eps\right) - W_v\left(\frac{1}{2}-2\eps\right), W_v\left(\frac{1}{2}+2\eps\right) - W_v\left(\frac{3}{4}-2\eps\right)\right)\geq \varphi_l(\eps).
\end{equation*}
 \end{Assumption}

Assumption~\ref{ass:quantile_fct_empirical_mean} is a pretty weak assumption since, intuitively, it requires that the distribution of the centered and variance one  real-valued random variables $\inr{\Sigma^{-1/2}(\overline{\tilde X}_i-\mu), v}$ have their $1/4$-quartiles and  medians constant far away as well as for their $3/4$-quartiles and medians, and this has to hold uniformly in all directions $v\in\cS_2^{d-1}$. For instance, in the Gaussian case, Assumption~\ref{ass:quantile_fct_empirical_mean} holds for $\varphi_u(\eps)=\Phi^{-1}(3/4)+c_0\eps$ and $\varphi_l(\eps)=\Phi^{-1}(3/4)-c_0\eps$ for some absolute constant $c_0$ and for all $0<\eps<1/10$  (where we recall that $\Phi:t\to \bP[g\leq t]$ where $g\sim\cN(0,1)$).

\begin{Proposition}\label{prop:MOMAD} There are absolute constants $c_0,c_1$ and $c_2$ such that the following holds. We assume that Assumption~\ref{ass:quantile_fct_empirical_mean} holds for some $0<\eps<1/8$ and constants $\varphi_l(\eps)$ and $\varphi_u(\eps)$. We assume that the adversarial contamination with $L_2$ inliers model from Assumption~\ref{assum:first} holds with a number of adversarial outliers $|\cO|\leq \eps K$. We assume that $K\geq c_0 \eps^{-2} d$. With probability at least $1-\exp(-c_1 \eps^2 K)$, for all $v\in\bR^d$,
\begin{equation*}
\varphi_l(\eps)\sqrt{\frac{K}{N}}\norm{\Sigma^{1/2}v}_2\leq MOMAD_K(v)  \leq \varphi_u(\eps)\sqrt{\frac{K}{N}}\norm{\Sigma^{1/2}v}_2.
\end{equation*} 
\end{Proposition}

Proposition~\ref{prop:MOMAD} shows that $MOMAD_K$ is equivalent to $v\to\sqrt{K/N}\norm{\Sigma^{1/2}v}_2$ up to the two constants $\varphi_u(\eps)$ and $\varphi_l(\eps)$. We will be interested in two situations regarding these constants. The first one is when their ratio is upper bounded by some absolute constant: there exists some absolute constant $c_0$ such that
\begin{equation}\label{eq:ratio_varphi_small}
\frac{\varphi_u(\eps)}{\varphi_l(\eps)}\leq c_0.
\end{equation}This condition will be enough to obtain robust optimal subgaussian bounds for $\hat \mu_{MOM,K}^{SDO}$ in the two following theorems. If condition \eqref{eq:ratio_varphi_small} holds we say that $MOMAD_K$ is isomorphic to $v\to\sqrt{K/N}\norm{\Sigma^{1/2}v}_2$. The second condition, that will be of interest to us is when we will estimate $\Sigma$ using $MOMAD_K$ in Section~\ref{sec:estimation_of_}, is when the two constants $\varphi_u(\eps)$ and $\varphi_l(\eps)$ can be made arbitrarily close to the same constant by taking $\eps$ small enough, that is when there exists some absolute constants $\phi_0$ and $c_1>0$ such that for all $0<\eps<\phi_0/c_1$, 
\begin{equation}\label{eq:varphi_cond}
 \varphi_l(\eps)=\phi_0-c_1\eps \mbox{ and } \varphi_u(\eps)=\phi_0+c_1\eps.
  \end{equation} In that case, we speak about an \textit{almost-isometric property} of $MOMAD_K$. The later condition is stronger than an isomorphic property but it allows to solve a higher order moment problem. In Section~\ref{sec:study_of_the_}, we provide several examples where these conditions hold as well as other properties of the family of cdfs $(H_v:v\in\cS_2^{d-1})$ even when there is not even a first moment.

We finish this section with an isomorphic result for $SDO_K$. The rate of convergence appears in this result: it is the level $r^*$ above which $SDO_K$ is isomorphic to $\nu\in\bR^d\to \norm{\Sigma^{-1/2}(\nu-\mu)}_2/\sqrt{K/N}$. One can define it as a solution to
\begin{equation}\label{eq:key_fixed_point_1}
C_0\left(\sqrt{\frac{d+1}{K}}+\sqrt{\frac{u}{K}}\right) + \sup_{\norm{v}_2=1}H_{N,K,v}(r^*) + \frac{|\cO|}{K}< \frac{1}{2}
\end{equation}where $u$ is a confidence parameter and $C_0$ a constant appearing in \eqref{eq:VC_concentration}.

  \begin{Proposition}\label{prop:iso_SDO}There are absolute constants $c_0,c_1$ and $c_2$ such that the following holds. We assume that Assumption~\ref{ass:quantile_fct_empirical_mean} holds for some $0<\eps<1/8$ and constants $\varphi_l(\eps)$ and $\varphi_u(\eps)$. 
We assume that the adversarial contamination with $L_2$ inliers model from Assumption~\ref{assum:first} holds with a number of adversarial outliers denoted by $|\cO|$. We assume that $|\cO|\leq \eps K$ and $K\geq c_0 \eps^{-2} d$. Let  $u>0$ and $r^*$ be such that \eqref{eq:key_fixed_point_1} holds. Then, with probability at least $1-\exp(-u)-\exp(-c_1 \eps^2 K)$, for all $\nu\in\bR^d$, if $\norm{\Sigma^{-1/2}(\nu-\mu)}_2\geq 2\sqrt{K/N}r^*$ then
\begin{equation*}
\frac{\norm{\Sigma^{-1/2}(\nu-\mu)}_2}{2\varphi_u(\eps)\sqrt{K/N}}  \leq SDO_K(\nu)\leq \frac{3\norm{\Sigma^{-1/2}(\nu-\mu)}_2}{2\varphi_l(\eps)\sqrt{K/N}} 
\end{equation*}and if $\norm{\Sigma^{-1/2}(\nu-\mu)}_2\leq 2\sqrt{K/N}r^*$ then $ SDO_K(\nu)\leq (3/\varphi_l(\eps))r^*$.
\end{Proposition}
Proposition~\ref{prop:iso_SDO} may be seen as a MOM version holding in the heavy-tailed case of the Proposition~\ref{prop:MAD} obtained in the Gaussian case. Such an extension from the Gaussian case to the $L_2$ heavy-tail case is made possible only thanks to the median-of-means principle and the use of the bucketed means instead of the data themselves. However, we will identify situations where condition~\eqref{eq:ratio_varphi_small} and \eqref{eq:key_fixed_point_1} with an optimal choice of rate $r^*$ (that is for the subgaussian rate \eqref{eq:subgaussian_rate_Sigma}) hold for $K=N$ even when a first moment do not exist. In that case, one can get a contamination price down to $|\cO|/N$ instead of the information theoretic lower bound in the general $L_2$ case given by $\sqrt{|\cO|/N}$ (see Section~\ref{sec:the_l2_case_via_behavior_in_zero}). We start with the general $L_2$ case and then we will consider an extra assumption that allows for such better bounds.

\subsection{The general $L_2$ case} 
\label{sec:the_L2_case_via_markov}
Unlike in Section~\ref{sec:the_gaussian_case} or Section~\ref{sec:the_l2_case_via_behavior_in_zero} below where, we demand that for all $v\in\cS_2^{d-1}$ and for values of $0<r<c_0$  the deviation function $H_{N,K,v}(r)$ is less than $1/2-c_1r$ here we simply use Markov inequality to control the function $H_{N,K,v}$ around $0$. The price we pay by using this approach is that we will not have anymore estimation results for the SDO MOM over $K$ blocks which hold for all deviation parameter $u$ up to $K$ but only for  $u\sim K$. The other price we pay here is for the adversarial contamination cost that will be of the order of $\sqrt{|\cO|/N}$ whereas as proved by Theorem~\ref{thm:SDO_L2_robust_subgauss} below it can be better up to $|\cO|/N$ (as in the Gaussian case from Theorem~\ref{thm:SDO_gauss_robust_subgauss}). We will be able to achieve this result thanks to an extra regularity assumption of the cdfs $H_v$ of all one-dimensional projections around $0$ (see Assumption~\ref{ass:cdf_empirical_close_0} below). But, for the moment, we do not grant this type of assumption in this section and obtain a general result under only the existence of a second moment as well as Assumption~\ref{ass:quantile_fct_empirical_mean}. Subgaussian rates can be derived out of this result when condition \eqref{eq:ratio_varphi_small} holds (we refer to Section~\ref{sec:study_of_the_} where this condition is studied). 
 
 In this section, the bound we use is simply the one deduced from Markov's inequality that is for all $r>0$ and $K\in[N]$: 
 \begin{equation}\label{eq:markov}
 H_{N,K,v}(r) = \bP\left[\frac{1}{\sqrt{N/K}}\sum_{i=1}^{N/K} \inr{\Sigma^{-1/2}(\tilde X_i-\mu), v}\geq r\right] \leq \frac{1}{1+r^2}.
 \end{equation}(Note that we used a slightly modification of Markov's inequality: if $Z$ is a centered variance one real-valued random variable then $\bP[Z\geq r]=\min_{a\in\bR}\bP[Z+a\geq r+a]\leq (1+r^2)^{-1}$). Our main result in the general $L_2$ setup will follow from this bound and a general result stated in Section~\ref{sec:proofs}. It is now stated in the following theorem.

\begin{Theorem}\label{thm:SDO_L2_2} There are absolute constants $c_0,c_1$ and $c_2$ such that the following holds. We assume that Assumption~\ref{ass:quantile_fct_empirical_mean} holds for some $0<\eps<1/8$.
We assume that the adversarial contamination with $L_2$ inliers model from Assumption~\ref{assum:first} holds with a number of adversarial outliers $|\cO|\leq c_0\eps K$. We assume that $c_0K\geq \eps^{-2}d$. With probability at least $1-2\exp(-c_0\eps^2 K)$,  
\begin{equation*}
\norm{\Sigma^{-1/2}(\hat \mu^{SDO}_{MOM,K} - \mu)}_2\leq \frac{4\varphi_u(\eps)}{\varphi_l(\eps)} \sqrt{\frac{K}{N}}.
\end{equation*}
\end{Theorem}

The rate of convergence in Theorem~\ref{thm:SDO_L2_2} can be written like the one in Theorem~\ref{thm:SDO_gauss_robust_subgauss} and Theorem~\ref{thm:SDO_L2_robust_subgauss} below where the three terms: complexity, deviation and price for adversarial corruption appear. Indeed, one should notice here that the deviation probability in Theorem~\ref{thm:SDO_L2_2} is fixed equal to $1-2\exp(-c_0\eps^2 K)$ because we had to take the deviation parameter $u$ equal to $K$ because of the approach based on Markov's inequality \eqref{eq:markov}. It is however, equivalent to replace $\sqrt{K/N}$ by $\sqrt{d/(\eps^2N)} + \sqrt{u/N}+\sqrt{|\cO|/(\eps N)}$ for $u=K$ since the two quantities are equivalent under the assumptions of Theorem~\ref{thm:SDO_L2_2}. In that case, one may recognize the complexity term $\sqrt{d/N}$, the deviation term $\sqrt{u/N}$ as well as the price for adversarial corruption $\sqrt{|\cO|/N}$. In particular, we see that the price we pay for the corruption is of the order of  $\sqrt{|\cO|/N}$ which is larger than the $|\cO|/N$ term in the Gaussian case from Theorem~\ref{thm:SDO_gauss_robust_subgauss} and it is the worst case of Theorem~\ref{thm:SDO_L2_robust_subgauss} below. Indeed, in  Theorem~\ref{thm:SDO_L2_2} we did not exploit any other property than the existence of a second moment whereas the two other two Theorems \ref{thm:SDO_gauss_robust_subgauss} and Theorem~\ref{thm:SDO_L2_robust_subgauss} exploit some regularity assumption around $0$ of the family of functions $H_{N,K,v}, v\in\cS_2^{d-1}$.


\paragraph{Adaptation to $K$ via Lepski's method.} It follows from Theorem~\ref{thm:SDO_L2_2} that $\hat \mu^{SDO}_{MOM,K}$ is an estimator which depends on the deviation parameter. Therefore, we need to construct an adaptive to $K$ version of this estimator to disentangle the estimator from the deviation parameter. The classical way to do it is via Lepski's method \cite{MR1091202,MR1147167}. Usually, the price we pay to make this approach work is some extra knowledge on $\Sigma$ such as its trace and operator norm. But here for the SDO type estimator we are using together with the isomorphic property of the $SDO_K$ we only need knowledge on $\varphi_u(\eps)$ and $\varphi_l(\eps)$. Let us now construct this adaptive scheme: the number of blocks is chosen via
\begin{equation}\label{eq:adapt_K_lepski}
\hat K = \min\left(K\in[N]: SDO_{k}(\hat \mu^{SDO}_{MOM,K} - \hat \mu^{SDO}_{MOM,k})\leq \max\left(\frac{9}{\varphi_l(\eps)}, \frac{6\varphi_u(\eps)}{\varphi_l^2(\eps)}\left(1+\sqrt{\frac{K}{k}}\right)\right), \forall k=N,\cdots, K\right).
\end{equation}

\begin{Theorem}\label{thm:SDO_L2_2_lepskii} There are absolute constants $c_0,c_1$ and $c_2$ such that the following holds. We assume that Assumption~\ref{ass:quantile_fct_empirical_mean} holds for some $0<\eps<1/8$ and all $K\in[N]$.
We assume that the adversarial contamination with $L_2$ inliers model from Assumption~\ref{assum:first} holds with a number of adversarial outliers denoted by $|\cO|$. Then, for all $K\geq \max(c_0\eps^{-2}d, c_0|\cO|)$ with probability at least $1-2\exp(-c_0\eps^2 K)$,  
\begin{equation*}
\norm{\Sigma^{-1/2}(\hat \mu^{SDO}_{MOM,\hat K} - \mu)}_2\leq \frac{28\varphi_u^2(\eps)}{\varphi_l^2(\eps)}\sqrt{\frac{K}{N}}.
\end{equation*}where $\hat K$ is the adaptive choice of number of blocks from \eqref{eq:adapt_K_lepski}. 
\end{Theorem}

\subsection{The  $L_2$ case under an extra regularity condition around $0$ of the $H_v$'s} 
\label{sec:the_l2_case_via_behavior_in_zero}
In this section, we obtain estimation bound for the MOM version of the SDO median in the adversarial corruption with heavy-tailed $L_2$ inliers model under an extra assumption on the regularity at $0$ of the family of functions $H_v, v\in\cS_2^{d-1}$ that is stated now.  

\begin{Assumption}\label{ass:cdf_empirical_close_0}
There exists some absolute constants $c_0,c_1>0$ and $c_2>0$ such that for all $v\in\cS_2^{d-1}$ and all $(2C_0/c_1)\sqrt{(d+1)/K}\leq r\leq c_0$
\begin{equation*}
H_{N,K,v}(r)= H_v(r) := \bP\left[\frac{1}{\sqrt{N/K}}\sum_{i=1}^{N/K} \inr{\Sigma^{-1/2}(\tilde X_i-\mu), v}\geq r\right] \leq \frac{1}{2}-c_2r.
\end{equation*}
\end{Assumption}
This assumption is about the behavior around the origin of the cdf of all one-dimensional projections of the random vectors $(N/K)^{-1/2}\sum_{i=1}^{N/K} \Sigma^{-1/2}(\tilde X_i-\mu)$ where the $\tilde X_i$ are the non-corrupted data. The term $\frac{1}{2}-c_2r$ in the bound above is the behavior of regular in $0$ cdfs such as in the Gaussian case (see Section~\ref{sec:study_of_the_} for more details and more examples).

Our main result in the adversarial corruption and  heavy-tailed $L_2$ model under Assumption~\ref{ass:cdf_empirical_close_0} is the following theorem. The proof may be found in Section~\ref{sec:proofs}.
\begin{Theorem}\label{thm:SDO_L2_robust_subgauss} There are absolute constants $c_0,c_1$ and $c_2$ such that the following holds. We assume that Assumption~\ref{ass:quantile_fct_empirical_mean} holds for some $0<\eps<1/4$ and that Assumption~\ref{ass:cdf_empirical_close_0} holds as well. We assume that the adversarial contamination with $L_2$ inliers model from Assumption~\ref{assum:first} holds with a number of adversarial outliers $|\cO|\leq c_0\eps K$. We assume that $c_0K\geq \eps^{-2}d$. For all $0<u\leq c_0\eps^2 K$, with probability at least $1-2\exp(-u)$, 
\begin{equation}\label{eq:thm_SDO_L2_robust_sub_main}
\norm{\Sigma^{-1/2}(\hat \mu^{SDO}_{MOM,K} - \mu)}_2\leq \frac{c_1\varphi_u(\eps)}{\varphi_l(\eps)} \left(\sqrt{\frac{d}{N}}+\sqrt{\frac{u}{N}}+ \frac{|\cO|}{\sqrt{NK}} \right).
\end{equation}
\end{Theorem}

We recover the optimal subgaussian rate \eqref{eq:subgaussian_rate_Sigma}  in Theorem~\ref{thm:SDO_L2_robust_subgauss} when for some $0<\eps<1/8$, condition~\eqref{eq:ratio_varphi_small} holds and $|\cO|\lesssim \sqrt{K d}$. The term $|\cO|/\sqrt{KN}$ appearing in the convergence rate of Theorem~\ref{thm:SDO_L2_robust_subgauss} is the price we pay for the adversarial contamination. It is between $\sqrt{|\cO|/N}$ when $K\sim|\cO|$ and $|\cO|/N$ when $K \sim N$. Usually when the inliers are only in $L_2$, the information theoretic lower bound is known to be of the order of $\sqrt{|\cO|/N}$  and not like $|\cO|/N$. We get a better rate in Theorem~\ref{thm:SDO_L2_robust_subgauss} thanks to Assumption~\ref{ass:cdf_empirical_close_0} which is  using in some more efficient way the regularity of the $H_v$ functions at $0$.

Unlike typical results in the MOM literature except for the one obtained in \cite{MS}, the deviation rate in Theorem~\ref{thm:SDO_L2_robust_subgauss} is  $1-2\exp(-u)$ for all $u\lesssim K$, in particular it does not have to depend on parameter $K$. As a consequence, the estimator $\hat \mu^{SDO}_{MOM,K}$ does not depend on the deviation parameter. Usually, results for MOM estimators constructed on $K$ blocks are given with probability at least $1-\exp(-c_0K)$ and then a Lepski's method is used to construct an adaptive to $K$ procedure. This is not the case here nor it is for the Gaussian case in Section~\ref{sec:the_gaussian_case} (as we did in the previous section). This is again because  Assumption~\ref{ass:cdf_empirical_close_0} is using more efficiently the behavior of $H_{N,K,v}$ around $0$.

\section{Estimation of $\Sigma$ under a $L_2$-moment assumption with MOMAD} 
\label{sec:estimation_of_}
In this section, we show that it is possible to estimate the covariance matrix $\Sigma$ using the MOMAD estimator. In particular, given that the isomorphic property of MOMAD hold under Assumption~\ref{ass:quantile_fct_empirical_mean} which does not require more moment than $L_2$ moment, we show that it is possible to estimate $\Sigma$ without requiring more moment than $2$ that is just under the assumption that $\Sigma$ exists. This differs from approaches based on the empirical covariance matrix where at best a $L_{2+\delta}$-moment assumption for some positive $\delta$ is granted for the estimation of the covariance matrix \cite{MR3217437,MR3476606,lu2020robust}. 

We show that for the estimation of $\Sigma$ via the MOMAD, the properties of $\varphi_l(\eps)$ and $\varphi_u(\eps)$ introduced in Assumption~\ref{ass:quantile_fct_empirical_mean} play a key role. Let us first have a look at these quantities in the Gaussian case. In that case, there are some absolute constants $\phi_0$ and $c_1>0$ such that for all $0<\eps<\phi_0/c_1$, 
\begin{equation}\label{eq:varphi_cond}
 \varphi_l(\eps)=\phi_0-c_1\eps \mbox{ and } \varphi_u(\eps)=\phi_0+c_1\eps
  \end{equation}where $\phi_0=\Phi^{-1}(3/4)$ (see Section~\ref{sec:study_of_the_} or the proof of Proposition~\ref{prop:MAD} for more details). This later result holds in the Gaussian case first because the two interquartile intervals have the same length: $\Phi^{-1}(0) - \Phi^{-1}(1/4)=\Phi^{-1}(3/4)-\Phi(0)=\phi_0$ and, second, because the Gaussian density function is uniformly lower bounded by an absolute positive  constant locally around the two $1/4$ and $3/4$ quartiles $\Phi^{-1}(1/4)$ and $\Phi^{-1}(3/4)$ as well as around the median $\Phi^{-1}(1/2)=0$. If this last condition were not true at some $q\in\{W(1/4), W(1/2), W(3/4)\}$ where $W=W_{N,K,v}$ for some direction $v\in\cS_2^{d-1}$  then there will be some plateau of the cdf $r\in\bR\to 1-H_{N,K,v}(r)$ starting at $q$ and thus there would be a constant factor gap between $W(\ell/4-2\eps)$ and $W(\ell/4+2\eps)$ for some $\ell\in\{1, 2, 3\}$. In that case, there would be some absolute constants $c_0>0$ and $c_1>0$ such that $|\varphi_l(\eps) - \varphi_l(\eps)|\geq c_0$ for all $0<\eps\leq c_1$. In particular, we would only have an isomorphic property for the MOMAD and thus its is not clear how to estimate $\Sigma$ using MOMAD at a better rate than a constant rate. Typical values of $\phi_0$ in \eqref{eq:varphi_cond} will be $\phi_0=W(1/4)-W(1/2)=W(1/2)-W(3/4)$. In particular, the interquartile interval lengths have to be equal in all directions $v\in\cS_2^{d-1}$; this will hold, in particular, under a spherical symmetry assumption of the $\Sigma^{-1/2}(\tilde X_i - \mu)$ (see Section~\ref{sec:study_of_the_} for more formal statement).

  \begin{Assumption}\label{ass:behabiour_varphi_fcts}For the same choice of $K$ as in Assumption~\ref{ass:quantile_fct_empirical_mean} where $\eps>0\to\varphi_l(\eps), \varphi_u(\eps)$ are defined, there are absolute constants $\phi_0$, $c_0>0$ and $c_1>0$ such that for all $v\in\cS_2^{d-1}$ and all $0<\eps<c_0$, $\varphi_l(\eps)=\phi_0-c_1\eps$ and $\varphi_u(\eps) = \phi_0+ c_1\eps$.
  \end{Assumption}

  Let us now turn to the construction of an estimator of the covariance matrix $\Sigma$ using MOMAD under Assumption~\ref{ass:behabiour_varphi_fcts} (as well as Assumption~\ref{ass:quantile_fct_empirical_mean}). Because of the constant factor $\phi_0$ in Assumption~\ref{ass:behabiour_varphi_fcts} we will provide an estimator of the \textit{scatter matrix} $\phi_0^2\Sigma$ (according to \cite{MR2238141}, a scatter matrix is any matrix proportional to the covariance matrix). 

  It follows from Proposition~\ref{prop:MOMAD} that $MOMAD_K$ is isomorphic to $v\in\bR^d\to\phi_0\sqrt{K/N}\norm{\Sigma^{1/2}v}_2$ and that under Assumption~\ref{ass:behabiour_varphi_fcts} it becomes an almost isometry, that is, with probability at least $1-\exp(-c_0\eps^2 K)$, for all $v\in\bR^d$, 
 \begin{equation}\label{eq:MOMAD_L2_cov_esti}
 \left|MOMAD_K(v) - \phi_0\sqrt{\frac{K}{N}}\norm{\Sigma^{1/2}v}_2\right|\leq c_1\eps\sqrt{\frac{K}{N}}\norm{\Sigma^{1/2}v}_2
 \end{equation}as long as $|\cO|\leq \eps K$ and $K\geq c_2\eps^{-2}d$. 
 In the Gaussian case and other spherical cases as in Section~\ref{sec:study_of_the_}, this almost isometric property holds for $K=N$ (and $MOMAD_N=MAD$) and any $0< \eps<1/4$: it follows from Proposition~\ref{prop:MAD} that with probability at least $1-\exp(-c_0\eps^2N)$, for all $v\in\bR^d$, 
\begin{equation}\label{eq:MAD_L2_cov_esti}
 \left|MAD(v) - \Phi^{-1}(3/4)\norm{\Sigma^{1/2}v}_2\right|\leq c_1 \eps\norm{\Sigma^{1/2}v}_2.
 \end{equation}

 We may use \eqref{eq:MOMAD_L2_cov_esti} to estimate directly the entries of $\Sigma$ following an idea from \cite{gnanadesikan1972robust}. Let $(e_j)_{j=1}^d$ denote the canonical basis of $\bR^d$. We have, for all $i,j\in[d]$,
 \begin{equation*}
  4\Sigma_{ij} = 4\inr{e_i, \Sigma e_j} = \norm{\Sigma^{1/2}(e_i+e_j)}_2^2 - \norm{\Sigma^{1/2}(e_i-e_j)}_2^2.
 \end{equation*}As a consequence, a natural estimator of $\phi_0^2\Sigma$ based on $MOMAD_K$ is the matrix $\hat \Sigma$ whose entries are defined for all $i,j\in[d]$ by
 \begin{equation*}
\hat\Sigma_{ij} = \frac{N}{4K}\left(MOMAD_K^2(e_i+e_j)- MOMAD_K^2(e_i-e_j)\right). 
 \end{equation*} Note that $\hat\Sigma$ is symmetric but it may not be a SDP. To overcome this issue, a projection method has been introduced in \cite{lu2020robust} which may also be used as well for $\hat \Sigma$. Our main statistical bound for $\hat \Sigma$ is the following.

 \begin{Proposition}\label{prop:estima_Sigma}Assume that Assumption~\ref{assum:first} holds. Let $K\in[N]$, $\varphi_l$ and $\varphi_u$ be such that Assumption~\ref{ass:quantile_fct_empirical_mean} and Assumption~\ref{ass:behabiour_varphi_fcts} hold. Then, for all $0<\eps<c_0$ such that $|\cO|\leq \eps K$ and $K\geq c_2 \eps^{-2}d$, with probability at least $1-\exp(-c_3\eps^2K)$, 
 \begin{equation*}
 \max_{i,j\in[d]}\left|\frac{\phi_0^2\Sigma_{ij} - \hat\Sigma_{ij}}{\Sigma_{ii} + \Sigma_{jj}}\right|\leq \sup_{\norm{u}_1=\norm{v}_1=1}\left|\frac{\inr{u,(\phi_0^2\Sigma - \hat \Sigma)v}}{\sum_{i}(|u_i|+|v_i|)\Sigma_{ii}}\right| \leq c_1\eps(c_1\eps + \phi_0)/2.
 \end{equation*}
 \end{Proposition}

 In particular, if one can choose $K=N$ so that Assumption~\ref{ass:quantile_fct_empirical_mean} and Assumption~\ref{ass:behabiour_varphi_fcts} hold -- for instance, in the Gaussian case or for other spherical variables as in Section~\ref{sec:study_of_the_} -- then the $MOMAD_N$ estimator becomes the classical MAD one and for $\eps^2 = c_2d/N$ we have that with probability at least $1-\exp(-c_4d)$,  
  \begin{equation*}
 \max_{i,j\in[d]}\left|\frac{\phi_0^2\Sigma_{ij} - \hat\Sigma_{ij}}{\Sigma_{ii} + \Sigma_{jj}}\right|\leq \sup_{\norm{u}_1=\norm{v}_1=1}\left|\frac{\inr{u,(\phi_0^2\Sigma - \hat \Sigma)v}}{\sum_{i}(|u_i|+|v_i|)\Sigma_{ii}}\right|\leq c_5\sqrt{\frac{d}{N}}.
 \end{equation*}as long as $|\cO|\leq c_6 d$.

\section{Study of the $H_{N,K,v}, v\in\cS_2^{d-1}$ functions} 
\label{sec:study_of_the_}
The functions $H_{N,K,v}, v\in\cS_{2}^{d-1}$ play a key role in our analysis. Their behavior in a neighborhood of their $1/4$ and $3/4$ quartiles and medians should be controlled so that Assumption~\ref{ass:quantile_fct_empirical_mean} may hold: they are driving the isomoprhic properties and almost isometric properties of the $MOMAD_K$ and $SDO_K$ functions and so of the statistical performance of the Stahel Donoho Median and Median of Means. Their behavior around $0$ also drives the improved rates obtained under Assumption~\ref{ass:cdf_empirical_close_0}. From our perspective, it is of the utmost importance to understand the behavior of these functions at these particular points. 

Let us first settle down the properties of the $H_{N,K,v}$ functions desirable for our analysis. We set $Z_i=\Sigma^{-1/2}(\tilde X_i-\mu)$ for all $i\in[N]$ so that the $Z_i$'s are independent centered isotropic vectors in $\bR^d$ and $n=N/K$.  We want to identify conditions on the distributions of the $Z_i$'s such that
\begin{itemize}
	\item for Assumption~\ref{ass:cdf_empirical_close_0}: there exists some absolute constants $c_0,c_1>0$ such that for all $v\in\cS_2^{d-1}$ and all $0<r<c_0$,
\begin{equation}\label{eq:recall_1}
H_{n,v}(r):=\bP\left[\frac{1}{\sqrt{n}}\sum_{i=1}^n \inr{Z_i, v}\geq r\right]\leq \frac{1}{2}-c_1 r.
 \end{equation}
 \item for Assumption~\ref{ass:quantile_fct_empirical_mean}: there exists some absolute constants $c_0>0$ and $0<\eps<1/8$ such that $\varphi_l(\eps)$ and  $\varphi_u(\eps)$ exist and are such that
 \begin{equation}\label{eq:recall_2}
 \frac{\varphi_u(\eps)}{\varphi_l(\eps)}\leq c_0 
 \end{equation}or there are absolute constants $\phi_0$ and $c_1>0$ such that for all $0<\eps<\phi_0/c_1$,
 \begin{equation}\label{eq:recall_3}
  \varphi_l(\eps)=\phi_0-c_1\eps \mbox{ and } \varphi_u(\eps)=\phi_0+c_1\eps
 \end{equation}which are respectively Condition~\ref{eq:ratio_varphi_small} (insuring an isomorphic property of $MOMAD_K$ and $SDO_K$ as well as optimal subgaussian rates for SD median and median of means) and Condition~\ref{eq:varphi_cond} (insuring almost isometric property of $MOMAD_K$ as well as estimation properties for $\hat \Sigma$ in Section~\ref{sec:estimation_of_}).
\end{itemize} Let us first study the Gaussian case which is our benchmark situation. We will then study other cases where the family of functions $H_{N,K,v}, v\in\cS_2^{d-1}$ satisfies these conditions.

\paragraph{The Gaussian case.} We recall that $\Phi:t\in\bR\to \bP[g\leq t]=\int_{-\infty}^t \phi(u)du$ where $\phi:u\in\bR\to (2\pi)^{-1/2}\exp(-u^2/2)$ is the Gaussian density function. We also denote $H_G:t\to 1-\Phi(t)$ and $W_G:p\in(0,1)\to H_G^{(-1)}(p)$ the inverse function of $H_G$ so that $W_G(p)=\Phi^{-1}(1-p)$. It follows from the mean value theorem that for all $t,\eps\in\bR$, $|H_G(t+\eps) - H_G(t)| \leq  \max(-\phi(t),-\phi(t+\eps))\eps$ so that around $0$ we have for all $c_0>0$ and $0<r<c_0$, $H_G(r)\leq 1/2 -\phi(c_0)r$. As a consequence, \eqref{eq:recall_1} holds in the Gaussian case for instance with $c_0=1$ and $c_1=\phi(1)$. Let us now look at the two other conditions in the Gaussian case. From the mean value theorem, we  have for all $p\in[1/2,1)$ and $\eps\geq0$ such that $p+\eps\in[1/2,1)$, $\eps/\phi(W_G(p)) \leq  W_G(p) - W_G(p+\eps)\leq  \eps/\phi(W_G(p+\eps))$ and for all $p\in(0,1/2]$ and $\eps\geq0$ such that $p-\eps\in(0,1/2]$, $\eps/\phi(W_G(p)) \leq  W_G(p-\eps)-W_G(p)\leq \eps/\phi(W_G(p+\eps))$. We conclude that there are absolute constants $c_0,c_1>0$ such that for all $0\leq \eps \leq c_0$, 
\begin{equation*}
\varphi_u(\eps) = \Phi^{-1}(3/4) + 2\eps\left(\frac{1}{\phi(\Phi^{-1}(3/4))} + \frac{1}{\phi(0)}\right) + c_1\eps^2
\end{equation*}and
\begin{equation*}
\varphi_l(\eps) = \Phi^{-1}(3/4) - 2\eps\left(\frac{1}{\phi(\Phi^{-1}(3/4))} + \frac{1}{\phi(0)}\right) - c_1\eps^2.
\end{equation*} So that both conditions \eqref{eq:recall_1} and \eqref{eq:recall_3} hold with $\phi_0 = \Phi^{-1}(3/4)$. In particular, we see that the values of the density function $\phi$ at the $1/4$ and $3/4$ quartiles (here we used that $\phi(\Phi^{-1}(3/4)) = \phi(\Phi^{-1}(1/4))$) and at the median $\phi(0)$ (here we used that $\Phi^{-1}(1/2)=0$) play a key role.

In the following, we identify situations where the $H_{N,K,v},v\in\cS_2^{d-1}$ functions and their pseudo inverses mimic the $H_G$ and $W_G$ functions from the Gaussian case. There are at least two reasons: the first one is that we are projecting random vectors leaving in $\bR^d$ onto one dimensional subspaces; the second reason is that we are averaging random variables having a second moment. We will explore these two observations in the two following paragraphs.

\paragraph{One dimensional projections and elliptically contoured distributions.} The fact that the $H_v$ functions deal only with one-dimensional marginals is making these functions likely to behave as in the Gaussian case  since one-dimensional projections of sufficiently spherically symmetric random vectors in $\bR^d$ are expected to behave like one-dimensional Gaussian variables and this phenomenon is even more accentuated when $d$ is large (this is one particular situation where large dimension $d$ may help in Statistics). Indeed, one may have in mind an observation -- sometimes attributed to H. Poincar{\'e} -- that the density function of the one-dimensional projection $\inr{\sqrt{d}U,e_1}$ -- where $\sqrt{d}U$ is uniformly distributed over $\sqrt{d}\cS_2^{d-1}$ and $(e_j)_{j=1}^d$ is the canonical basis of $\bR^d$ -- converges to the density of a $\cN(0,1)$ when $d\to \infty$ (see page~16 in \cite{MR2814399} or Chapter~4 in \cite{MR1335228}). One may also have in mind that there are directions such as $v=(1/\sqrt{d}, \ldots,1/\sqrt{d})$ which are mixing the coordinates of $\Sigma^{-1/2}(\tilde X_1-\mu)$ when projected onto $v$ and therefore may have the tendency to mimic a standard Gaussian variable because of the CLT. 
Note that all these observations hold for $N=K$ that is even for $n=1$: because of the one-dimensional projections we may not even have to average the $Z_i$'s to mimic the Gaussian case. Therefore, Theorem~\ref{thm:SDO_general} can be extended beyond the Gaussian case when this phenomenon occurs. 

Let us now consider an example of elliptically contoured distributions where this happens to be true. Our aim is to show that Condition~\eqref{eq:ratio_varphi_small} and Assumption~\ref{ass:cdf_empirical_close_0} (and so Theorem~\ref{thm:SDO_L2_robust_subgauss}) may hold for $K=N$ (i.e. $n=1$) even when the $\tilde X_i$'s do not have a first moment. 

We assume that the $\tilde X_i$'s are i.i.d. and that $\Sigma^{-1/2}(\tilde X_1-\mu)$ has a spherically symmetric distribution; in that case, $\tilde X_1-\mu$ is sometimes said to have an \textit{elliptically contoured  distribution}. Then, there exists a non-negative random variable $R$ such that  $\Sigma^{-1/2}(\tilde X_1-\mu)$ is distributed according to $R U$ where $U$ is uniformly distributed on $\cS_2^{d-1}$ and is independent of $R$ (see Chapter~4 in \cite{MR1335228}). In that case, all the $\inr{\Sigma^{-1/2}(\tilde X_1-\mu), v}$ for $v\in\cS_2^{d-1}$ have the same distribution as  $\inr{\Sigma^{-1/2}(\tilde X_1-\mu), e_1}$ (where $(e_j)_{j=1}^d$ is the canonical basis of $\bR^d$) which is distributed according to $R\inr{U,e_1}$. Now, using that $\inr{U,e_1}$ is absolutely continuous w.r.t. the Lebesgue measure with density function given by, when $d\geq2$, 
\begin{equation*}
t\in\bR\to C_d (1-t^2)^{\frac{d-3}{2}}I(|t|\leq 1) \mbox{ where } C_d=\left(\int_{-1}^1 (1-t^2)^{\frac{d-3}{2}}dt\right)^{-1} =\frac{2\Gamma(d/2)}{\Gamma((d-1)/2)\sqrt{\pi}}
\end{equation*} and $\Gamma$ is the Gamma function, we can deduce that (even for $K=N$), $H_v$ is independent of $v\in\cS_2^{d-1}$ and is such that for all $r\geq 0$, $H_v(-r)=1-H_v(r)$ and 
\begin{equation*}
H_v(r)=H(r) := C_d\int_{0}^1 \bP[R\geq r/x] \left(1-x^2\right)^{\frac{d-3}{2}}dx.
\end{equation*}In particular, we recover that $H(0)=1/2$ since $R\geq0$ a.s.. Let us now consider a simple example for the distribution of $R$. In that example, $R$ takes values $r_1<r_2<\cdots$ such that $\alpha_j = \bP[R=r_j]$ for all $j\in\bN^*$ so that for all $q>0$, $\bE R^q = \sum_j r_j^q\alpha_j$ which may be infinite even for $q=1$ (that is when there is not even a first moment). For this example, we have for all $r\geq0$,
\begin{equation*}
H(r)=C_d\sum_{j=1}^\infty \alpha_j \int_{r/r_j}^1 (1-x^2)^{\frac{d-3}{2}}dx I(r\leq r_j).
\end{equation*}In particular, $H$ is differentiable and $R\inr{U,e_1}$ is absolutely continuous w.r.t. the Lebesgue measure with a density function given by 
\begin{equation*}
f:r\in\bR\to - H^\prime(r)=C_d\sum_{j=1}^\infty \frac{\alpha_j}{r_j} \left[1-\left(\frac{r}{r_j}\right)^2\right]^{\frac{d-3}{2}} I(r\leq r_j).
\end{equation*}In particular, for $r_\infty = \lim_{j\to \infty} r_j$, $H$ is strictly decreasing on $[0, r_\infty)$ from $H(0)=1/2$ to $H(r_\infty)=0$ and beyond $r_\infty$ it is constant equal to $0$. Therefore, for all $v\in\cS_2^{d-1}$, the generalized inverse $W_v$ of $H_v$ is independent of $v$ and it is the  inverse of $H$: for all $p\in(0,1/2]$ there is a unique element $W(p) (=W_v(p))$ in $[0,r_\infty)$ such that $H(W(p))=p$ and $W(1-p)=-W(p)$. 

Now, let us choose $r_j=2^j C_d$ and $\alpha_j=2^{-j}$ for all $j\in\bN$. We also assume $d\geq4$ to make the presentation simpler (the cases $d=1,2,3$ can be treated separately). In that case, $\bE R=+\infty$ and so the mean and covariance matrix do not exist. Nevertheless, one may still assume that there exists $\mu\in\bR^d$ and $\Sigma\in\bR^{d\times d}$ definite positive such that $\Sigma^{-1/2}(\tilde X_1 - \mu)$ is spherically symmetric (without having $\mu$ to be a mean vector and $\Sigma$ to be a covariance matrix). Then, Theorem~\ref{thm:SDO_L2_robust_subgauss} still applies.

Let us first check Condition~\eqref{eq:recall_1}. We have $H(0)=1/2$ and for all $0\leq r\leq C_d/\sqrt{d-3}$,
\begin{equation}\label{eq:cond_around_0_checked}
 f(r)\geq \sum_{j=1}^\infty \frac{1}{2^{2j}} \left[1-{\frac{d-3}{2C_d^2}}\left(\frac{r}{2^j}\right)^2\right]\geq \frac{1}{3}.
 \end{equation} Moreover, we see that $\sqrt{d}\leq C_d\leq 6\sqrt{d}$, hence, \eqref{eq:cond_around_0_checked} holds for all $0\leq r\leq 1$. Which is according to the means value theorem enough to show that Condition~\eqref{eq:recall_1} holds (see \eqref{eq:mean_value_theo} below for more details).

Let us now check conditions \eqref{eq:recall_2} and \eqref{eq:recall_3}. It follows from Proposition~\ref{prop:sym_density} below that, it is enough to lower bound the density function $f$ in a neighborhood of $p$ for $p\in\{W(1/4), W(1/2), W(3/4)\}$ and that $W(1/4)-W(3/4)$ is an absolute constant. But, given that $W(1/2)=0$ and \eqref{eq:cond_around_0_checked} holds, that $f$ is symmetric about $0$ and that $W(1/4) = -W(3/4)$, we only have to checked that $f(q)\geq c_0$ for all $q\in[W(1/4)-2\eps, W(1/4)+2\eps]$ for some $0<\epsilon<1/8$ and an absolute constant $c_0$ and that $W(1/4)$ is an absolute constant. We first have to find $W(1/4)$ which is the unique solution $r$ such that $H(r)=1/4$. We see that $f$ is symmetric unimodal with maximal value at $0$ given by $f(0) =4/3$ and we showed that $f(r)\geq 1/3$ for all $0\leq r\leq 1$ in \eqref{eq:cond_around_0_checked}. Therefore, $H(1/8)\geq 1/3$ and $H(1-1/10)\leq 3/15< 1/4$, hence, $W(1/4)\in[1/8, 1-1/10]$. It follows from \eqref{eq:cond_around_0_checked} that $f(q)\geq 1/3$ for all $q\in[W(1/4)-1/10, W(1/4)+1/10]$. We conclude that both conditions \eqref{eq:recall_2} and \eqref{eq:recall_3} hold thanks to Proposition~\ref{prop:sym_density} below.

For this example, one can take $\varphi_u(\eps) = W(1/4)-(4/3)\eps$ and $ \varphi_l(\eps) = W(1/4)+(4/3)\eps$ for all $0<\eps<1/16$. In that case, $MOMAD_K$ is an almost isometry and we can state a result like Theorem~\ref{thm:SDO_gauss_robust_subgauss} where $\Phi^{-1}(3/4)$ is replaced by $W(1/4)$ and $\mu$ and $\Sigma$ are not anymore the mean and covariance matrix since they do not exist but 'location' and 'scale' parameters defined such that $\Sigma^{-1/2}(\tilde X_i - \mu)$ are spherically symmetric. As a consequence, the phenomenon underlying the Gaussian case from Section~\ref{sec:the_gaussian_case} has nothing to do with concentration but it is more about elliptical symmetry.   



\paragraph{Gaussian approximation.} 
In cases where there is some lack of spherical symmetry of $\Sigma^{-1/2}(\tilde X_1-\mu)$ one may study the $H_v$ functions for a smaller number $K$ of blocks so that $n=N/K$ may be large enough to see some averaging effect. In that case and because Gaussian variables satisfy all the properties we need, it is tempting to use a Gaussian approximation result such as a Berry-Esseen bound (see \cite{MR1353441,MR3289373,MR1841329}) to approximate the $H_v$ functions by $1-\Phi$  for $n=N/K$ large enough. This strategy has been used several times in Minsker and co-authors works on Median-of-means and Catoni's type of estimators (see for instance \cite{MS,minsker2018uniform}).


For instance, when for all $v\in\cS_2^{d-1}$,  $\inr{\Sigma^{-1/2}(\tilde X_i-\mu),v},i\in[n]$ are (independent, centered and variance one) real-valued random variables in $L_{2+\delta}$ such that  $\norm{\inr{\Sigma^{-1/2}(\tilde X_i-\mu),v}}_{2+\delta}\leq\kappa$ (uniformly in $v\in\cS_2^{d-1}$) for some $\delta>0$ then, it follows from  Theorem~5.7 in \cite{MR1353441} that there is an absolute constant $c_0>0$ such that for all $v\in\cS_2^{d-1}$ and all $r\in\bR$, 
\begin{equation}\label{eq:BE_L_2_plus_delta}
\left|H_{N,K,v}(r)  - \bP[g\geq r]\right|\leq \frac{c_1\kappa^{2+\delta}}{n^{\delta/2}}:=c_n
\end{equation}It follows that for all $p\in(0,1)$ and $\eps\in\bR$ satisfying $p+\eps\in(0,1)$ that
\begin{equation*}
\Phi^{-1}\left(1-p-\eps-c_n\right)\leq W(p+\eps)\leq \Phi^{-1}\left(1-p-\eps + c_n\right).
\end{equation*}In particular, for all $0<\eps<1/16$, if $n$ is large enough so that $c_n\leq \eps$ then one can take $\varphi_u(\eps)=\Phi^{-1}(3/4)-c_0\eps$ and $\varphi_l(\eps)=\Phi^{-1}(3/4) - c_0\eps$. So that the ratio $\varphi_u(\eps)/\varphi_l(\eps)$ is constant; in that case, the $MOMAD_K$ and $SDO_k$ are isomorphism (see Proposition~\ref{prop:MOMAD}) and we recover a subgaussian rate in Theorem~\ref{thm:SDO_L2_robust_subgauss}. 

However, a Gaussian approximation result such as the one in \eqref{eq:BE_L_2_plus_delta} is not enough for Assumption~\ref{ass:cdf_empirical_close_0}. Indeed, it follows from \eqref{eq:BE_L_2_plus_delta} that for all $0\leq r\leq c_0$, $H_v(r) \leq H_G(r) + c_n\leq 1/2-c_1r+c_n$ for some absolute constants $c_0>0$ and $c_1>0$. It appears that our analysis used to prove Theorem~\ref{thm:SDO_L2_robust_subgauss} does not stand this extra error term $c_n$ compare with Assumption~\ref{ass:cdf_empirical_close_0}. Gaussian approximation does not help in this case: indeed Assumption~\ref{ass:cdf_empirical_close_0} is more about the existence of a uniform lower bound around $0$ of the density functions of the one-dimensional projections $\inr{n^{-1/2}\sum_i Z_i, v}$ as we are considering now.

\paragraph{Beyond the Gaussian behavior.} In the later two paragraphs, we identified situations where the $n^{-1/2}\sum_{i=1}^n \inr{Z_i, v}$ for $v\in\cS_2^{d-1}$ behave like Gaussian variables. We saw that this may be the case because we are considering one-dimensional projections of $d$-dimensional vectors and/or we are taking empirical means over $n$ variables. But properties we are looking for the $H_{n,v},v\in\cS_2^{d-1}$ functions (see \eqref{eq:recall_1}, \eqref{eq:recall_2} and \eqref{eq:recall_3}) are all dealing only with their behavior around $3$ (or $4$ when the median is not $0$) points. So that only the behavior of these functions at these points play a role and there is no need to mimic the Gaussian case for all values of $r$ in $\bR$. We now state a general result going in this direction. In particular, we recover the conditions from \cite{maronna1995behavior} and \cite{MR2051003} recalled in the Introduction section.

Let us assume that the $n^{-1/2}\sum_{i=1}^n \inr{Z_i, v}$ for $v\in\cS_2^{d-1}$ are absolutely continuous w.r.t. the Lebesgue measure with a density function denoted by $f_v$. By the mean value theorem, we have for all $r\geq0$, all $p\in(0,1)$ and $\eps\geq0$ such that $p+\eps\in(0,1)$,
\begin{equation}\label{eq:mean_value_theo}
H_v(r)\leq H_v(0) - \min_{0\leq t\leq r}f_v(t) r \mbox{ and } \frac{\eps}{\max_{q\in[p,p+\eps]}f_v(W_v(q))} \leq W_v(p) - W_v(p+\eps)\leq  \frac{\eps}{\min_{q\in[p,p+\eps]}f_v(W_v(q))}.
\end{equation}In particular, the values of the density functions $f_v,v\in\cS_2^{d-1}$ at $0,W(1/4), W(1/2)$ and $W(3/4)$ drives the quality of inequalities from \eqref{eq:mean_value_theo} and, as noted in previous works on the Stahel-Donoho outlyingness function, are enough to insure all the conditions we need on $H_v$ and $W_v$ recalled in \eqref{eq:recall_1}, \eqref{eq:recall_2} and \eqref{eq:recall_3}.

\begin{Proposition}\label{prop:sym_density} Let $K\in[N]$ be such that $N/K\in\bN$. We assume that the original non-corrupted data $\tilde X_i,i\in[N]$ are independent and that there exists $\mu\in\bR^d$ and $\Sigma\in\bR^{d\times d}$ definite positive so that for all $v\in\cS_2^{d-1}$,  $\sqrt{K/N}\sum_{i=1}^{N/K}\inr{\Sigma^{-1/2}(\tilde X_i-\mu),v}$ are absolutely continuous real valued random variables with a density denoted by $f_v$. 

If there exists $0<\eps<1/8$ and $c_0>0$ such that for all $v\in\cS_2^{d-1}$, all $p\in\{W_v(1/4), W_v(1/2), W_v(3/4)\}$ and all $q\in[p-2\eps, p+2\eps]$, $f_v(q)\geq c_0$ then for $\cI_v=\max\left(W_v(1/4)-W_v(1/2), W_v(1/2)-W_v(3/4)\right)$
\begin{equation*}
\varphi_u(\eps) = \max_{v\in\cS_2^{d-1}}\left(\cI_v + 4\eps/c_0\right) \mbox{ and } \varphi_l(\eps) = \min_{v\in\cS_2^{d-1}}\left(\cI_v - 4\eps/c_0\right).
\end{equation*}
 We also have
\begin{equation*}
\frac{\varphi_u(\eps)}{\varphi_l(\eps)}\leq \frac{\max_{v\in\cS_2^{d-1}}\cI_v}{\min_{v\in\cS_2^{d-1}}\cI_v}\left(1+\frac{16\eps}{c_0\min_{v\in\cS_2^{d-1}}\cI_v}\right)
 \end{equation*}when $4\eps\leq c_0\min_v \cI_v$. Moreover, if $(c_0/4)\max_v \cI_v<1/8$ and  $\min_v\cI_v\geq c_1$, for some absolute constant $c_1>0$,  then  condition \eqref{eq:recall_2} holds (and so we recover the optimal subgaussian rates in Theorem~\ref{thm:SDO_L2_2} and Theorem~\ref{thm:SDO_L2_2_lepskii}) and if for all $v\in\cS_2^{d-1}$, $\cI_v:=\phi_0$ then condition~\eqref{eq:recall_3} holds and so does Proposition~\ref{prop:estima_Sigma}.

If for all $v\in\cS_2^{d-1}$, $H_v(0)\leq 1/2$ and there are absolute constants $c_0>0$ and $c_1>0$ so that for all $0<r<c_0$, $f_v(v)\geq c_1$ then Assumption~\ref{ass:cdf_empirical_close_0} holds (that is \eqref{eq:recall_3} holds) and so does Theorem~\ref{thm:SDO_L2_robust_subgauss}.
\end{Proposition}


Note that in Proposition~\ref{prop:sym_density}, $\mu$ and $\Sigma$ do not have to be the mean and covariance matrix of the $\tilde X_i$'s. In that case, $\mu$ and $\Sigma$ are sometimes called location and scale and so Theorem~\ref{thm:SDO_L2_robust_subgauss} still apply for the robust to adversarial contamination and heavy-tail estimation of location, even in situations where there is not even a first moment.

Proposition~\ref{prop:sym_density} gives an alternative to Gaussian approximation which does not, in general, allow to check Assumption~\ref{ass:cdf_empirical_close_0} because  of the residual terms in Esseen or Berry-Esseen type inequalities. The assumptions in Proposition~\ref{prop:sym_density} are all granting that the density functions $f_v$ are locally lower bounded around the 'critical' $1/4$ and $3/4$ quartiles and medians.  They are natural assumptions that already appeared in several studies of estimators based on the SDO. In Proposition~\ref{prop:sym_density} we show that by using the median-of-means principle these assumptions are dealing with the density functions on the bucketed means and not the data themselves. However, Proposition~\ref{prop:sym_density} may also be applied in the $K=N$ case as for  elliptically contoured distributions.

 \section{Conclusion} 
\label{sec:conclusion}
We showed that it is possible to estimate a mean vector in $\bR^d$ w.r.t. the metric $\norm{\Sigma^{-1/2}\cdot}_2$ even though $\Sigma$ is unknown, the data set is corrupted by an adversary and the data are heavy-tailed. The rate obtained are the (deviation) minmax one in the ideal i.i.d. Gaussian case. The estimator used to achieve this rate is a deepest point with respect a median-of-means version of the Stahel-Donoho outlyingness functional. When the data are spherical enough there is no need to bucket the data and then the estimator is using the classical  Stahel-Donoho outlyingness. Our analysis shows that the two cases can be handled using the same methodology and that the family of cdfs $(H_{N,K,v}:v\in\cS_2^{d-1})$ plays a key role in this analysis, in particular, their behavior around $0$, the median and the $1/4$ and $3/4$ quartiles. 

In this work, we have not deal with several research opportunities opened by the SDO. We now list some of them that may be considered in future works. a) It may look possible to use the isomorphic properties of the $MOMAD_K$ and $SDO_K$ to study the Stahel-Donoho estimator (SDE) or a median-of-means version of the SDE defined as 
\begin{equation}\label{eq:MOM_SDE}
\tilde \mu^{SDE}_{MOM,K} = \frac{\sum_{k=1}^K \hat w_k \bar X_k}{\sum_{k=1}^K \hat w_k}
\end{equation}where $(\hat w_k)_{k=1}^K$ are non-negative weights such that $\hat w_k$ depends on the outlyingness of the $k$-th bucketed mean $\bar X_k$. For instance,
\begin{equation}
 \hat w_k = \left\{\begin{array}{cc}
1 & \mbox{ if } SDO_K(\bar X_k)\leq \hat \alpha_K\\
0 & \mbox{ otherwise.} 
\end{array} \right. \mbox{ where } \hat \alpha_k = \Med(SDO_K(\bar X_k)).
 \end{equation} b) Similarly, the isomorphic or almost-isometry properties of $MOMAD_K$ and $SDO_K$ may also be used to study the properties of a MOM version of the SDE of the covariance matrix:
\begin{equation*}
\hat \Sigma = \frac{2}{K}\sum_{k=1}^K \hat w_k (\bar X_k- \tilde \mu^{SDE}_{MOM,K})(\bar X_k- \tilde \mu^{SDE}_{MOM,K})^\top.
\end{equation*} c) From a computational point of view, it is still an open question to construct an approximate solution to the SDO. The original or MOM version of the Stahel-Donoho median could be approximated via a robust gradient descent algorithms such as the one introduced in \cite{Bartlett19,Jules_Guillaume_1,lei2020fast} with some extra normalization step required by the MAD denominator. We expect this algorithm to be more efficient than the classical weighted SDE because we expect to do only $\log d$ iterations to achieve a subgaussian estimator using a robust gradient descent algorithm whereas the SDE would require to approximate the $K$ depths $SDO_K(\bar X_k), k\in[K]$  and should therefore require more computational time (note that, in practice the SDE has been reported to be more efficient than the deepest data that is the data $\bar X_k$ with the smallest $SDO_K(\bar X_k)$ but the SDE was not compared with an approximate solution of $\hat \mu_{MOM,K}^{SDO}$).


\section{Proofs} 
\label{sec:proofs}
In this section, we provide some proofs of all the results from the preceding sections. The only complexity measure we are using in this work is the Vapnik and Chervonenkis (VC) dimension \cite{MR3408730,MR1719582} of a class $\cF$ of Boolean functions, i.e. of functions from $\bR^d$ to $\{0,1\}$ in our case. We recall that $VC(\cF)$ is the maximal integer $n$ such that there exists $x_1, \ldots, x_n\in\bR^d$ for which the set $\{(f(x_1), \cdots, f(x_n)):f\in\cF)\}$ is of maximal cardinality that is $2^n$. The only VC-dimension we will use is the one of the set of all indicators of half affine spaces in $\bR^d$: $VC(\{x\in\bR^d\to I(\inr{\cdot, v}\geq r):v\in\bR^d, r\in\bR\}) = d+1$ (see Example~2.6.1 in \cite{MR1385671}). The main technical tool (see Chapter~3 in \cite{Kol11}) we will be using is the following one: let $Y_1, \ldots, Y_n$ be independent random vectors in $\bR^d$, there exists an absolute constant $C_0$ such that for all $u>0$, with probability at least $1-\exp(-u)$,
\begin{equation}\label{eq:VC_concentration}
\sup_{f\in\cF}\left(\frac{1}{n}\sum_{i=1}^n f(Y_i)-\bE f(Y_i)\right)\leq C_0\left(\sqrt{\frac{VC(\cF)}{n}} + \sqrt{\frac{u}{n}}\right).
\end{equation}

 We recall that for all $v\in\cS_2^{d-1}$, $K\in[N]$ and $r>0$,
\begin{equation*}
H_{N,K,v}(r) = \bP\left[\frac{1}{\sqrt{N/K}}\sum_{i=1}^{N/K} \inr{\Sigma^{-1/2}(\tilde X_i-\mu), v}\geq r\right].
\end{equation*}The rate of convergence we will obtained is any $r^*$ satisfying 
\begin{equation}\label{eq:key_fixed_point}
C_0\left(\sqrt{\frac{d+1}{K}}+\sqrt{\frac{u}{K}}\right) + \sup_{\norm{v}_2=1}H_{N,K,v}(r^*) + \frac{|\cO|}{K}< \frac{1}{2}
\end{equation}where $C_0$ is the constant from \eqref{eq:VC_concentration} and for some choice of $K$ and $u$ specified in each result depending on the set of assumptions.

\subsection{Proof of Proposition~\ref{prop:MOMAD} and \ref{prop:MAD} (first part): isomorphic property of MOMAD} 
\label{sub:proof_of_the_props}
We first prove Proposition~\ref{prop:MOMAD} -- the proof of Proposition~\ref{prop:MAD} is a straighforward application of Proposition~\ref{prop:MOMAD}.

\textbf{Proof of Proposition~\ref{prop:MOMAD}.}
We first observe that by renormalization, it is enough to show that for all $v\in\cS_2^{d-1}$,
\begin{equation}\label{eq:prop_MOMAD_main_1}
\varphi_l(\eps)\leq \Med(|\inr{\Sigma^{-1/2}(\bar X_k-\mu), v} - \Med(\inr{\Sigma^{-1/2}(\bar X_k-\mu),v})|) \leq \varphi_u(\eps).
\end{equation} Moreover, for all $i\in[N], \Sigma^{-1/2}(\tilde X_i-\mu)$ has mean zero and covariance $I_d$. Hence, without loss of generality we assume that $\mu=0$ and $\Sigma=I_d$.

The strategy we are using to prove \eqref{eq:prop_MOMAD_main_1} is the following one. Let $K$ real numbers $a_1, \ldots, a_K$ be given and denote by $a_{(1)}\leq \cdots \leq a_{(K)}$ the non-decreasing rearrangement of the $(a_k)_k$ (this is the rearrangement of the $a_k$'s and not of their absolute values).  To prove a result like $\varphi_l(\eps)\leq \Med(|a_k - \Med(a_k)|)\leq \varphi_u(\eps)$, it is enough to show that $\varphi_l(\eps)\leq a_{(3(K+1)/4)} - a_{((K+1)/2)}\leq \varphi_u(\eps)$ and $\varphi_l(\eps)\leq a_{((K+1)/2)} - a_{((K+1)/4)}\leq \varphi_u(\eps)$. As a consequence, to prove a result like \eqref{eq:prop_MOMAD_main_1}, we should study the rearrangement (the two quartiles and the median) of the $\inr{\bar X_k,v}, k\in[K]$ uniformly over all $v\in\cS_2^{d-1}$.  But, $|\cO|$ elements among the $X_i$'s come from the adversary and we do not have any control on their behavior. We therefore have to consider the worst possible case which is when $|\cO|$ bucketed means $\bar X_k$ are corrupted by one outliers from $\{X_i:i\in\cO\}$. However, one may check that if we change $|\cO|$ points in a set $\{a_k:k\in[K]\}$ to get a new set $\{A_k:k\in[K]\}$ then   $\varphi_l(\eps)\leq a_{(3(K+1)/4)} - a_{((K+1)/2)}\leq \varphi_u(\eps)$ will be true if we show that $\varphi_l(\eps)\leq A_{(3(K+1)/4-|\cO|)} - A_{((K+1)/2+|\cO|)}$  and $A_{(3(K+1)/4+|\cO|)} - A_{((K+1)/2-|\cO|)}\leq \varphi_u(\eps)$ -- and a similar observation holds for the other $(1/4)$-quartile. We will therefore first study the rearrangement of the original (i.e. non corrupted) bucketed means (later denoted by $\overline{\tilde X}_k, k\in[K]$) projected on all one dimensional directions uniformly over these directions to deduce the result from \eqref{eq:prop_MOMAD_main_1} on the corrupted bucketed means $\bar X_k$.

We denote by $\overline{\tilde X}_k,k\in[K]$ the bucketed means of the original (non corrupted) dataset, i.e. $\overline{\tilde X}_k = (1/|B_k|)\sum_{i\in B_k}\tilde X_i$ for $k\in[K]$. To prove \eqref{eq:prop_MOMAD_main_1} we first study the rearrangements of vectors $(\inr{\overline{\tilde X}_k, v})_{k\in[K]}$ uniformly over all $v\in\cS_2^{d-1}$. We will then deal with the adversarial corruption to get \eqref{eq:prop_MOMAD_main_1}.
 
We introduce the following supremum of empirical process:
\begin{equation*}
Z=\sup_{\ell\in[K-1]}\sup_{\norm{v}_2=1}\left|\frac{1}{K}\sum_{k=1}^K I\left(\inr{\overline{\tilde X}_k,v}\geq \frac{W_v(\ell/K)}{\sqrt{N/K}}\right) - \bP\left[\inr{\overline{\tilde X}_k,v}\geq \frac{W_v(\ell/K)}{\sqrt{N/K}}\right]\right|
\end{equation*}where $W_v$ has been defined in Definition~\ref{def:H_and_W}. It follows from \eqref{eq:VC_concentration} that for all $u>0$, with probability at least $1-\exp(-u)$, $Z\leq C_0\left(\sqrt{(d+1)/K} + \sqrt{u/K}\right)$ (note that even though the function $W_v$ depends on $v$, the boolean function $x\to I(\inr{x,v}\geq W_v(\ell/N))$ is still the indicator of a affine half-space of $\bR^d$ for all $v\in\bR^d$ and all $\ell\in[K-1]$ and thus the VC dimension of the set of Boolean functions $\{x\to I(\inr{x, v}\geq W_v(\ell/K)):v\in\bR^d, \ell\in[K-1]\}$ is less or equal to $d+1$).  As a consequence, for some choice of  $0<\eps<1/8$ such that Assumption~\ref{ass:quantile_fct_empirical_mean} holds, if $K\geq (2C_0)(d+1)\eps^{-2}$ then with probability at least $1-\exp(-\eps^2 K/(2C_0)^2)$, $Z\leq \eps$. Let us denote by $\Omega_\eps$ the event onto which $Z\leq \eps$; we proved that $\bP[\Omega_\eps]\geq 1-\exp(-\eps^2 K/(2C_0)^2)$.

Let us place ourselves on the event $\Omega_\eps$ up to the end of the proof. Since for all $v\in\cS_2^{d-1}$,
\begin{equation*}
\bP\left[\inr{\overline{\tilde X}_k,v}\geq \frac{W_v(\ell/K)}{\sqrt{N/K}}\right] = H_v(W_v(\ell/K))=\ell/K,
\end{equation*}(by left continuity of $H_v$ we have $H_v(W_v(p))=p$ for all $p\in(0,1)$), we have for all $\ell\in[K]$ and $v\in\cS_2^{d-1}$, that 
\begin{equation}\label{eq:cardinal_esti_inliers}
\left|\left\{k\in[K]: \inr{\overline{\tilde X}_k,v}\geq \frac{W_v(\ell/K)}{\sqrt{N/K}}\right\}\right|\in\left[\ell- \eps K, \ell+\eps K\right].
\end{equation}This last result on the uniform in $v\in\cS_2^{d-1}$ rearrangement of $(\inr{\overline{\tilde X}_k,v})_k$ will be used to get the desire result on the rearrangement for $(\inr{\overline{X}_k,v})_k$ (uniformly in $v$). To go from the $\overline{\tilde X}_k$'s to the $\overline{X}_k$'s we now have to deal with the adversarial corruption.

Since, there are $|\cO|$ original data that may have been modified by the adversary, in the worse case $|\cO|$ bucketed means $\overline{\tilde X}_k$ may be considered as corrupted and so, from the above cardinality estimation result \eqref{eq:cardinal_esti_inliers}, we may only certify (on $\Omega_\eps$) that 
\begin{equation*}
\left|\left\{k\in[K]: \inr{\bar X_k,v}\geq \frac{W_v(\ell/K)}{\sqrt{N/K}}\right\}\right|\in\left[\ell- \eps K-|\cO|, \ell+\eps K+|\cO|\right]\subset \left[\ell- 2\eps K, \ell+2\eps K\right]
\end{equation*}on the $K$ bucketed means $\bar X_k$ constructed from the adversarialy corrupted dataset $\{X_i:i\in[N]\}$. We used here the assumption that $|\cO|\leq \eps K$. If follows from the later result that if we denote by $q^{1/4}_{K,v}$ the $1/4$ quartile of vector $(\inr{\bar X_k,v}:k\in[K])$, by $q^{1/2}_{K,v}$ its median and by $q^{3/4}_{K,v}$ its $3/4$ quartile then,
\begin{equation*}
 \sqrt{\frac{K}{N}}W_v\left(\frac{3}{4}+2\eps\right)\leq q^{1/4}_{K,v} \leq \sqrt{\frac{K}{N}}W_v\left(\frac{3}{4}-2\eps\right) ; \sqrt{\frac{K}{N}}W_v\left(\frac{1}{2}+2\eps\right)\leq q^{1/2}_{K,v} \leq  \sqrt{\frac{K}{N}}W_v\left(\frac{1}{2}-2\eps\right)
\end{equation*}and 
\begin{equation*}
\sqrt{\frac{K}{N}}W_v\left(\frac{1}{4}+2\eps\right)\leq q^{3/4}_{K,v}\leq \sqrt{\frac{K}{N}}W_v\left(\frac{1}{4}-2\eps\right).
\end{equation*}It follows from these inequalities that on the event $\Omega_\eps$, we have for all $v\in\cS_2^{d_1}$,
\begin{equation*}
\Med(|\inr{\bar X_k, v} - \Med(\inr{\bar X_k,v})|) \leq \sqrt{\frac{K}{N}}\max\left(W_v\left(\frac{1}{4}-2\eps\right) - W_v\left(\frac{1}{2}+2\eps\right), W_v\left(\frac{1}{2}-2\eps\right) - W_v\left(\frac{3}{4}+2\eps\right)\right)
\end{equation*}and 
\begin{equation*}
\Med(|\inr{\bar X_k, v} - \Med(\inr{\bar X_k,v})|) \geq \sqrt{\frac{K}{N}}\min\left(W_v\left(\frac{1}{4}+2\eps\right) - W_v\left(\frac{1}{2}-2\eps\right), W_v\left(\frac{1}{2}+2\eps\right) - W_v\left(\frac{3}{4}-2\eps\right)\right).
\end{equation*}The result follows from the definition of $\varphi_l(\eps)$ and $\varphi_u(\eps)$ in Assumption~\ref{ass:quantile_fct_empirical_mean}.
\endproof

\subsection{Proof of Proposition~\ref{prop:iso_SDO} and \ref{prop:MAD} (second part): isomorphic property of $SDO_K$.} 
\label{sub:isomorphic_property_of_}
The proof of Proposition~\ref{prop:iso_SDO} and \ref{prop:MAD} (second part) relies on the next result.

\begin{Proposition}\label{prop:numerator_SDO}
We assume that the adversarial contamination with $L_2$ inliers model from Assumption~\ref{assum:first} holds with a number of adversarial outliers denoted by $|\cO|$. Let $K\in[N]$,  $u>0$ and $r^*$ be such that \eqref{eq:key_fixed_point} holds. Then, with probability at least $1-\exp(-u)$,
\begin{equation*}
\sup_{v\in\cS_2^{d-1}}|\Med(\inr{\Sigma^{-1/2}(\bar X_k-\mu),v)})|\leq \sqrt{\frac{K}{N}}r^*. 
\end{equation*}
\end{Proposition}
\textbf{Proof of Proposition~\ref{prop:numerator_SDO}.} Denote by $\cK=\{k:B_k\cap\cO=\emptyset\}$ the set of indices of non-corrupted blocks of data. It follows from \eqref{eq:VC_concentration} and the definition of $r^*$ that with probability at least $1-\exp(-u)$, for all $v\in\cS_2^{d-1}$,
\begin{align*}
&\frac{1}{K}\sum_{k=1}^K I\left(\inr{\Sigma^{-1/2}(\bar X_k-\mu),v)}\geq \frac{r^*}{\sqrt{N/K}}\right)\\ 
&= \frac{1}{K}\sum_{k\in\cK}  I\left(\inr{\Sigma^{-1/2}(\overline{\tilde X}_k-\mu),v)}\geq  \frac{r^*}{\sqrt{N/K}}\right) + \frac{1}{K}\sum_{k\in \cK^c} I\left(\inr{\Sigma^{-1/2}(\bar X_k-\mu),v)}\geq  \frac{r^*}{\sqrt{N/K}}\right)\\
&\leq \frac{1}{K}\sum_{k=1}^K I\left(\inr{\Sigma^{-1/2}(\overline{\tilde X}_k-\mu),v)}\geq \frac{r^*}{\sqrt{N/K}}\right) + \frac{|\cO|}{K}\\
&\leq \sup_{\norm{v}_2=1} \left(\frac{1}{K}\sum_{k=1}^K I\left(\inr{\Sigma^{-1/2}(\overline{\tilde X}_k-\mu),v)}\geq \frac{r^*}{\sqrt{N/K}}\right) - P\left(\inr{\Sigma^{-1/2}(\overline{\tilde X}_k-\mu),v)}\geq \frac{r^*}{\sqrt{N/K}}\right)\right)\\ 
&+ P\left(\inr{\Sigma^{-1/2}(\overline{\tilde X}_1-\mu),v)}\geq \frac{r^*}{\sqrt{N/K}}\right)+ \frac{|\cO|}{K}\\
&\leq C_0 \left(\sqrt{\frac{d+1}{K}}+\sqrt{\frac{u}{K}}\right) + H_{N,K,v}(r^*) + \frac{|\cO|}{K}< \frac{1}{2}.
\end{align*}As a consequence,  with probability at least $1-\exp(-u)$, for all $v\in\cS_2^{d-1}$,
\begin{equation*}
\sum_{k=1}^K I\left(\inr{\Sigma^{-1/2}(\bar X_k-\mu),v)}\geq \frac{r^*}{\sqrt{N/K}}\right) < \frac{K}{2}
\end{equation*}and so
\begin{equation}\label{eq:SDO_L2_3}
\sup_{v\in\cS_2^{d-1}}|\Med(\inr{\Sigma^{-1/2}(\bar X_k-\mu),v)})|\leq \sqrt{\frac{K}{N}}r^*. 
\end{equation}

\begin{Remark}
It is also possible to consider a ''directional version'' of Proposition~\ref{prop:numerator_SDO} if one defines a ''directional version'' of $r^*$, that is for all directions $v\in\cS_2^{d-1}$, define $r^*_v>0$ satisfying
\begin{equation*}
C_0\left(\sqrt{\frac{d+1}{K}}+\sqrt{\frac{u}{K}}\right) + H_{N,K,v}(r^*_v) + \frac{|\cO|}{K}< \frac{1}{2}.
\end{equation*}Then, under the same conditions as in Proposition~\ref{prop:numerator_SDO}, we have with probability at least $1-\exp(-u)$,
\begin{equation*}
\sup_{v\in\cS_2^{d-1}}\frac{|\Med(\inr{\Sigma^{-1/2}(\bar X_k-\mu),v)})|}{r^*_v}\leq \sqrt{\frac{K}{N}}. 
\end{equation*}Hence, Proposition~\ref{prop:numerator_SDO} holds as well for $r^*=\sup_{\norm{v}_2=1}r^*_v$. Note that for most of the $v\in\cS_2^{d-1}$ the values of $r^*_v$ is expected to be much smaller than $r^*$. For instance, for vectors $v$ well-spread, we expect them to have a strong ''mixing'' power (see for instance ''super-Gaussian directions'' in \cite{MR3645123} or \cite{MR3136463,MR2520120}).
\end{Remark}

\paragraph{Proof of Proposition~\ref{prop:iso_SDO} and \ref{prop:MAD}.}  It follows from Proposition~\ref{prop:numerator_SDO} and Proposition~\ref{prop:MOMAD} that, with probability at least $1-\exp(-u)-\exp(-c_1 \eps^2 K)$, for all $v\in\cS_2^{d-1}$,
\begin{equation*}
|\Med(\inr{\Sigma^{-1/2}(\bar X_k-\mu),v)})|\leq \sqrt{\frac{K}{N}}r^* 
\end{equation*}and
\begin{equation*}
\varphi_l(\eps)\sqrt{\frac{K}{N}}\norm{\Sigma^{1/2}v}_2\leq MOMAD_K(v)  \leq \varphi_u(\eps)\sqrt{\frac{K}{N}}\norm{\Sigma^{1/2}v}_2.
\end{equation*}We denote by $\Omega_0$ the event onto which the last two properties hold. On the even $\Omega_0$, for all $\nu\in\bR^d$, we have
\begin{align*}
SDO_K(\nu) &= \sup_{v\in\bR^d}\frac{|\Med(\inr{\bar X_k-\nu,v)})|}{MOMAD_K(v)}\leq \sup_{v\in\bR^d}\frac{|\Med(\inr{\bar X_k-\nu,v)})|}{\varphi_l(\eps)\sqrt{K/N}\norm{\Sigma^{1/2}v}_2}=\sup_{v\in\cS_2^{d-1}}\frac{|\Med(\inr{\Sigma^{-1/2}(\bar X_k-\nu),v)})|}{\varphi_l(\eps)\sqrt{K/N}}\\
&\leq  \sup_{v\in\bR^d}\frac{|\Med(\inr{\Sigma^{-1/2}(\bar X_k-\mu),v})| + |\inr{\Sigma^{-1/2}(\nu-\mu),v}|}{\varphi_l(\eps)\sqrt{K/N}}\leq \sup_{v\in\bR^d}\frac{\sqrt{K/N}r^*+ |\inr{\Sigma^{-1/2}(\nu-\mu),v}|}{\varphi_l(\eps)\sqrt{K/N}}\\
&\leq  \left\{
\begin{array}{cc}
\frac{3\norm{\Sigma^{-1/2}(\nu-\mu)}_2}{2\varphi_l(\eps)\sqrt{K/N}}  & \mbox{ if } \norm{\Sigma^{-1/2}(\nu-\mu)}_2\geq 2\sqrt{K/N}r^*\\
3r^*/\varphi_l(\eps) & \mbox{ otherwise}
\end{array}\right.
\end{align*}and when $\norm{\Sigma^{-1/2}(\nu-\mu)}_2\geq 2\sqrt{K/N}r^*$, we have
\begin{equation*}
SDO_K(\nu)\geq \sup_{v\in\cS_2^{d-1}}\frac{|\inr{\Sigma^{-1/2}(\nu-\mu),v}| - |\Med(\inr{\Sigma^{-1/2}(\bar X_k-\mu),v})|}{\varphi_u(\eps)\sqrt{K/N}}\geq \frac{\norm{\Sigma^{-1/2}(\nu-\mu)}_2}{2\varphi_u(\eps)\sqrt{K/N}}. 
\end{equation*}
\endproof

\subsection{Proof of the statistical bounds} 
\label{sub:proof_of_the_statistical_bounds}
\paragraph{Proof of Proposition~\ref{prop:MAD}.}  Proposition~\ref{prop:MAD} is a corollary of Proposition~\ref{prop:MOMAD} for $K=N$. For this choice of $K$, there are $N$ blocks, each containing only one data and so $MOMAD_N(v)=MAD(v)$ for all $v\in\bR^d$. The only thing that remains to be checked is the validity of Assumption~\ref{ass:quantile_fct_empirical_mean} in the Gaussian case and the dependency of the $\varphi_l(\eps)$ and $\varphi_u(\eps)$  in terms of $\eps$.

When the original data $\tilde X_i,i\in[N]$ are $N$ i.i.d. Gaussian vectors $G_1, \ldots, G_N$ with mean $\mu$ and covariance matrix $\Sigma$ then for all $K\in[N]$, $(1/\sqrt{N/K})\sum_{i=1}^{N/K} \Sigma^{-1/2}(\tilde X_i-\mu)$ is a standard Gaussian vector in $\bR^d$. Therefore the  $H:=H_{N,K,v}$ function from Assumption~\ref{ass:quantile_fct_empirical_mean} is equal to the function $x\in\bR\to 1-\Phi(x)$ where  $\Phi:x\in\bR\to\bP[g\leq x]$ is the cdf of a standard Gaussian variable $g\sim\cN(0,1)$ in $\bR$. This holds for all $N,K$ and $v\in\cS_2^{d-1}$, that is $H_{N,K,v}$ is independent of $N,K$ and $v\in\cS_2^{d-1}$. Since $W:=W_{N,K,v}$ is the generalized inverse of $H$, in the Gaussian case, we obtain that $W(p)=\Phi^{-1}(1-p)$ for all $p\in(0,1)$. It follows from Lemma~5.2 in \cite{MR1353441} that there exists some absolute constant $C_1>0$ such that
\begin{equation*}
\min\left(W\left(\frac{1}{4}+2\eps\right) - W\left(\frac{1}{2}-2\eps\right), W\left(\frac{1}{2}+2\eps\right) - W\left(\frac{3}{4}-2\eps\right)\right)\geq \Phi^{-1}(3/4) - C_1 \eps:=\varphi_l(\eps)
\end{equation*}and 
\begin{equation*}
\max\left(W\left(\frac{1}{4}-2\eps\right) - W\left(\frac{1}{2}+2\eps\right), W\left(\frac{1}{2}-2\eps\right) - W\left(\frac{3}{4}+2\eps\right)\right)\leq \Phi^{-1}(1/4)\leq \Phi^{-1}(3/4) + C_1 \eps:=\varphi_u(\eps).
\end{equation*}As a consequence, Assumption~\ref{ass:quantile_fct_empirical_mean} holds in the Gaussian case for all $0<\eps<\Phi^{-1}(3/4)/C_1$ with $\varphi_l(\eps) = \Phi^{-1}(3/4) - C_1 \eps$ and  $\varphi_u(\eps) = \Phi^{-1}(3/4) + C_1 \eps$.
\endproof

\paragraph{Proofs of theorems~\ref{thm:SDO_gauss_robust_subgauss}, \ref{thm:SDO_L2_2} and \ref{thm:SDO_L2_robust_subgauss}} 
Theorems~\ref{thm:SDO_gauss_robust_subgauss}, \ref{thm:SDO_L2_2} and \ref{thm:SDO_L2_robust_subgauss} are corollaries of a general result that we are stating now.
\begin{Theorem}\label{thm:SDO_general} There are absolute constants $c_0,c_1$ and $c_2$ such that the following holds. We assume that Assumption~\ref{ass:quantile_fct_empirical_mean} holds for some $0<\eps<1/4$ and constants $\varphi_l(\eps)$ and $\varphi_u(\eps)$. We assume that the adversarial contamination with $L_2$ inliers model from Assumption~\ref{assum:first} holds with a number of adversarial outliers denoted by $|\cO|$. Let $K\geq \max(\eps^{-1}|\cO|,c_0\eps^{-2}d)$,  $0<u<c_0\eps^2 K$ and $r^*$ be such that \eqref{eq:key_fixed_point} holds. Then, with probability at least $1-2\exp(-u)$, 
\begin{equation*}
\norm{\Sigma^{-1/2}(\hat \mu^{SDO}_{MOM,K} - \mu)}_2\leq \frac{2\varphi_u(\eps)}{\varphi_l(\eps)} \sqrt{\frac{K}{N}}r^*.
\end{equation*}
\end{Theorem}

\paragraph{Proof of Theorem~\ref{thm:SDO_gauss_robust_subgauss}.}
There exists an absolute constant $c_0$ such that for all $0\leq r \leq c_0, \bP[g\geq r]\leq 1/2-2r$ where $g\sim\cN(0,1)$. Moreover, for all $K\in[N],v\in\cS_2^{d-1}$ and $r>0$, we have $H_{N,K,v}(r) = \bP[g\geq r]$. As a consequence, one can chose $r^*$, $u$ and $K$ such that
\begin{equation*}
r^* = C_0\left(\sqrt{\frac{d+1}{K}}+\sqrt{\frac{u}{K}}\right)+ \frac{|\cO|}{K} 
\end{equation*}as long as this later quantity is less or equal to $c_0$. Finally, we apply Theorem~\ref{thm:SDO_general} for $K=N$ and the result follows since $\hat \mu^{SDO}_{MOM,N} =\hat \mu^{SDO}$.
\endproof

\paragraph{Proof of Theorem~\ref{thm:SDO_L2_2}.} It follows from Markov's inequality \eqref{eq:markov} that we can chose $u$, $r^*$ and $K$ such that
\begin{equation*}
C_0\left(\sqrt{\frac{d+1}{K}}+\sqrt{\frac{u}{K}}\right) + \frac{1}{1+(r^*)^2} + \frac{|\cO|}{K}<\frac{1}{2}
\end{equation*}for instance by taking $r^*=2$, $K\geq 4|\cO|$, $K> 16C_0^2(d+1)$ and $K>16 C_0 u$. Note however, that because $r^*$ is constant, the convergence rate is proportional to $\sqrt{K/N}$, in particular it does not depend on $u$. Hence there is no interest to consider values of $u$ smaller than $K$ (up to constant).  We therefore apply Theorem~\ref{thm:SDO_general} for this choice of $K$, $u=c_0\eps^2K$ and $r^*=2$.
\endproof

\paragraph{Proof of Theorem~\ref{thm:SDO_L2_robust_subgauss}.}
Thanks to Assumption~\ref{ass:cdf_empirical_close_0}, there exists  absolute constants $c_0>0$ and $c_1>0$ such that for all $v\in\cS_2^{d-1}$ and $(2C_0/c_1)\sqrt{(d+1)/K}\leq r \leq c_0, H_{N,K,v}\leq 1/2-c_1r$. As a consequence, one can chose $r^*$, $u$ and $K$ such that
\begin{equation*}
r^* = \frac{2}{c_1}\left(C_0\left(\sqrt{\frac{d+1}{K}}+\sqrt{\frac{u}{K}}\right)+ \frac{|\cO|}{K} \right)
\end{equation*}as long as this later quantity is less or equal to $c_0$. Finally, we apply Theorem~\ref{thm:SDO_general} for this choice of $K$, $u$ and $r^*$.
\endproof

\paragraph{Proof of Theorem~\ref{thm:SDO_general}.} We first note that a proof of Theorem~\ref{thm:SDO_general} may follow from the isomorphic property of $SDO_K$ from Proposition~\ref{prop:iso_SDO}. However, it is possible to improve constants by using the following strategy.  

Let us place ourselves on the intersection of the two events where the results of both Proposition~\ref{prop:MOMAD} and Proposition~\ref{prop:numerator_SDO} hold.  We set $f:v\in\bR^d\to \Med(\inr{\bar X_k,v})$. Since $f$ is symmetric we have
\begin{align*}\label{eq:SDO_main_1}
\notag&\norm{\Sigma^{-1/2}(\hat \mu^{SDO}_{MOM,K} - \mu)}_2 =  \sup_{\norm{v}_2=1}\inr{\Sigma^{-1/2}(\hat \mu^{SDO}_{MOM,K} - \mu), v} = \sup_{v\in\bR^d} \inr{\hat \mu^{SDO}_{MOM,K} - \mu, \frac{v}{\norm{\Sigma^{1/2}v}_2}}\\
\notag & =\sup_{v\in\bR^d} \frac{\inr{\hat \mu^{SDO}_{MOM,K},v} - f(v) + f(v) -\inr{\mu,v}}{MOMAD_K(v)} \frac{MOMAD_K(v)}{\norm{\Sigma^{1/2}v}_2}\\
&\leq \left(\sup_{v\in\bR^d} \frac{\inr{\hat \mu^{SDO}_{MOM,K},v} - f(v) }{MOMAD_K(v)} + \sup_{v\in\bR^d} \frac{f(v) -\inr{\mu,v}}{MOMAD_K(v)}\right)\sup_{v\in\bR^d}\frac{MOMAD_K(v)}{\norm{\Sigma^{1/2}v}_2}\\
&\leq \left(SDO_K(\hat \mu^{SDO}_{MOM,K}) + SDO_K(\mu)\right)\sup_{v\in\bR^d}\frac{MOMAD_K(v)}{\norm{\Sigma^{1/2}v}_2}  \leq 2 SDO_K(\mu) \sup_{v\in\bR^d}\frac{MOMAD_K(v)}{\norm{\Sigma^{1/2}v}_2}.
\end{align*}where we used that $SDO_K(\hat \mu^{SDO}_{MOM,K})\leq SDO_K(\mu)$ by definition of $\hat \mu^{SDO}_{MOM,K}$.

We know how to control $\sup_{v\in\bR^d} MOMAD_K(v)/\norm{\Sigma^{1/2}v}_2$ by $\sqrt{K/N}\varphi_u(\eps)$ using Proposition~\ref{prop:MOMAD}. It remains to control the term $SDO_K(\mu)$. We have
\begin{align*}
\notag SDO_K(\mu) &= \sup_{v\in\bR^d}\frac{|\inr{\mu, v}-\Med(\inr{\bar X_k, v})|}{\Med(|\inr{\bar X_k, v} - \Med(\inr{\bar X_k,v})|)} = \sup_{v\in\bR^d} \frac{|\Med(\inr{\mu-\bar X_k,v})|}{\norm{\Sigma^{1/2} v}_2} \frac{\norm{\Sigma^{1/2}v}_2}{MOMAD_K(v)}\\
&\leq \sup_{\norm{v}_2=1} |\Med(\inr{\Sigma^{-1/2}(\bar X_k-\mu),v)})| \sup_{v\in\bR^d}\frac{\norm{\Sigma^{1/2}v}_2}{MOMAD_K(v)}.
\end{align*}The term $\sup_{v\in\bR^d}\norm{\Sigma^{1/2}v}_2 / MOMAD_K(v)$ is smaller than  $\sqrt{N/K}/\varphi_l(\eps)$ thanks to Proposition~\ref{prop:MOMAD}. Finally, to finish the proof, we  upper bound the term $\sup_{\norm{v}_2=1} |\Med(\inr{\Sigma^{-1/2}(\bar X_k-\mu),v)})|$ by $\sqrt{K/N}r^*$ thanks to Proposition~\ref{prop:numerator_SDO}.
\endproof

\paragraph{Proof of Theorem~\ref{thm:SDO_L2_2_lepskii}} 
\label{par:proof_of_theorem_Lepski}
For all $k\in[N]$, we set $\hat \mu_k = \hat \mu^{SDO}_{MOM,k}$ we denote by $\Omega_k$ the event onto which
\begin{equation*}
\norm{\Sigma^{-1/2}(\hat \mu_k - \mu)}_2\leq \frac{4\varphi_u(\eps)}{\varphi_l(\eps)} \sqrt{\frac{k}{N}}
\end{equation*}and, for all $\nu\in\bR^d$, if $\norm{\Sigma^{-1/2}(\nu-\mu)}_2\geq 6\sqrt{k/N}$ then
\begin{equation*}
\frac{\norm{\Sigma^{-1/2}(\nu-\mu)}_2}{2\varphi_u(\eps)\sqrt{k/N}}  \leq SDO_k(\nu)\leq \frac{3\norm{\Sigma^{-1/2}(\nu-\mu)}_2}{2\varphi_l(\eps)\sqrt{k/N}} 
\end{equation*}and if $\norm{\Sigma^{-1/2}(\nu-\mu)}_2\leq 6\sqrt{k/N}$ then
\begin{equation*}
 SDO_k(\nu)\leq \frac{9}{\varphi_l(\eps)}.
\end{equation*}
It follows from Proposition~\ref{prop:iso_SDO} for $r^*=3$ and $u=K/(16C_0^2)$ and Theorem~\ref{thm:SDO_L2_2} that $\bP[\Omega_k]\geq 1-3\exp(-c_1\eps^2 k)$ when $k\geq \max(|\cO|/\eps, c_0d/\eps^2)$. 

Let $K\geq \max(|\cO|/\eps, c_0d/\eps^2)$. On the even $\cap_{k=K}^N \Omega_k$, we have for all $K\leq k\leq N$,
\begin{align*}
&SDO_k(\hat \mu_K - \hat \mu_k)\leq \max\left(\frac{9}{\varphi_l(\eps)},\frac{3\norm{\Sigma^{-1/2}(\hat \mu_K - \hat \mu_k)}_2}{2\varphi_l(\eps)\sqrt{k/N}}\right)\leq\max\left(\frac{9}{\varphi_l(\eps)},\frac{6\varphi_u(\eps)}{\varphi_l^2(\eps)}\left(1+\sqrt{\frac{K}{k}}\right) \right)
\end{align*}and so, by definition of $\hat K$, we have $\hat K\leq K$. We also have by definition of $\hat K$ and because $\hat K\leq K$ that
\begin{equation*}
SDO_K(\hat \mu_{\hat K} - \hat \mu_K)\leq \max\left(\frac{9}{\varphi_l(\eps)},\frac{6\varphi_u(\eps)}{\varphi_l^2(\eps)}\left(1+\sqrt{\frac{\hat K}{K}}\right) \right)\leq \frac{12\varphi_u(\eps)}{\varphi_l^2(\eps)}.
\end{equation*}We conclude that either $\norm{\Sigma^{-1/2}(\hat\mu_{\hat K}-\hat \mu_K)}_2\leq 6\sqrt{K/N}$ and so  
\begin{equation*}
\norm{\Sigma^{-1/2}(\hat\mu_{\hat K}-\mu)}_2\leq \norm{\Sigma^{-1/2}(\hat\mu_{\hat K}-\hat \mu_K)}_2 + \norm{\Sigma^{-1/2}(\mu-\hat \mu_K)}_2\leq \left(6+\frac{4\varphi_u(\eps)}{\varphi_l(\eps)}\right)\sqrt{\frac{K}{N}}
\end{equation*}or $\norm{\Sigma^{-1/2}(\hat\mu_{\hat K}-\hat \mu_K)}_2\geq 6\sqrt{K/N}$ and so
\begin{align*}
\norm{\Sigma^{-1/2}(\hat \mu_{\hat K}-\mu)}_2  &\leq  \norm{\Sigma^{-1/2}(\hat \mu_{\hat K}-\hat\mu_K)}_2 + \norm{\Sigma^{-1/2}(\hat \mu_{K}-\mu)}_2\\
&\leq SDO_K(\hat \mu_{\hat K} - \hat \mu_K)2\varphi_u(\eps)\sqrt{\frac{K}{N}} + \frac{4\varphi_u(\eps)}{\varphi_l(\eps)} \sqrt{\frac{K}{N}}\leq \frac{28\varphi_u^2(\eps)}{\varphi_l^2(\eps)}\sqrt{\frac{K}{N}}.
\end{align*}
\endproof

\paragraph{Proof of Proposition~\ref{prop:estima_Sigma}.} We have for all $i,j\in[d]$, $\left|\phi_0^2\Sigma_{ij} - \hat\Sigma_{ij}\right|\leq c_1\eps(c_1\eps+\phi_0) \left(\Sigma_{ii}+\Sigma_{ij}\right) $ because, it follows from \eqref{eq:MOMAD_L2_cov_esti} that for all $v\in\bR^d$,
 \begin{align*}
 \left| MOMAD_K^2(v) - \phi_0^2\frac{K}{N}\norm{\Sigma^{1/2}v}_2^2\right| &= \left|MOMAD_K(v) - \phi_0 \sqrt{\frac{K}{N}}\norm{\Sigma^{1/2}v}_2\right|\left(MOMAD_K(v) + \phi_0 \sqrt{\frac{K}{N}}\norm{\Sigma^{1/2}v}_2\right)\\
 &\leq c_1 \eps \frac{K}{N}\norm{\Sigma^{1/2}v}_2^2(c_1\eps +\phi_0).
 \end{align*}

 Next, we have for all $u,v\in\bR^d$ such that $\norm{u}_1=\norm{v}_1=1$
 \begin{align*}
 &|\inr{u, (\phi_0^2 \Sigma - \hat \Sigma)v}| \\
 &= \frac{N}{4K} \left|\sum_{i,j} u_i v_j\left(\phi_0^2 \frac{K}{N}\norm{\Sigma^{1/2}(e_i+e_j)}_2^2 - MOMAD_K^2(e_i+e_j) + \phi_0^2 \frac{K}{N}\norm{\Sigma^{1/2}(e_i-e_j)}_2^2 - MOMAD_K^2(e_i-e_j) \right) \right|\\
 &\leq \frac{c_1 \eps(c_1\eps+\phi_0)}{4}\sum_{i,j} |u_i| |v_j| \left(\norm{\Sigma^{1/2}(e_i+e_j)}_2^2 + \norm{\Sigma^{1/2}(e_i-e_j)}_2^2\right) = \frac{c_1 \eps(c_1\eps+\phi_0)}{2}\sum_{i,j} |u_i| |v_j| \left(\Sigma_{ii} + \Sigma_{jj}\right).
 \end{align*}
\endproof

\begin{footnotesize}
\bibliographystyle{plain}
\bibliography{biblio}

\begin{thebibliography}{10}

\bibitem{MR1688610}
Noga Alon, Yossi Matias, and Mario Szegedy.
\newblock The space complexity of approximating the frequency moments.
\newblock {\em J. Comput. System Sci.}, 58(1, part 2):137--147, 1999.
\newblock Twenty-eighth Annual ACM Symposium on the Theory of Computing
  (Philadelphia, PA, 1996).

\bibitem{barany2020application}
Imre B{\'a}r{\'a}ny and Nabil~H Mustafa.
\newblock An application of the universality theorem for tverberg partitions to
  data depth and hitting convex sets.
\newblock {\em Computational Geometry}, page 101649, 2020.

\bibitem{MR1335228}
W\l~odzimierz Bryc.
\newblock {\em The normal distribution}, volume 100 of {\em Lecture Notes in
  Statistics}.
\newblock Springer-Verlag, New York, 1995.
\newblock Characterizations with applications.

\bibitem{MR3124669}
S\'{e}bastien Bubeck, Nicol\`o Cesa-Bianchi, and G\'{a}bor Lugosi.
\newblock Bandits with heavy tail.
\newblock {\em IEEE Trans. Inform. Theory}, 59(11):7711--7717, 2013.

\bibitem{MR3476606}
T.~Tony Cai, Weidong Liu, and Harrison~H. Zhou.
\newblock Estimating sparse precision matrix: optimal rates of convergence and
  adaptive estimation.
\newblock {\em Ann. Statist.}, 44(2):455--488, 2016.

\bibitem{MR3052407}
Olivier Catoni.
\newblock Challenging the empirical mean and empirical variance: a deviation
  study.
\newblock {\em Ann. Inst. Henri Poincar\'{e} Probab. Stat.}, 48(4):1148--1185,
  2012.

\bibitem{MR1841329}
Louis H.~Y. Chen and Qi-Man Shao.
\newblock A non-uniform {B}erry-{E}sseen bound via {S}tein's method.
\newblock {\em Probab. Theory Related Fields}, 120(2):236--254, 2001.

\bibitem{MR3845006}
Mengjie Chen, Chao Gao, and Zhao Ren.
\newblock Robust covariance and scatter matrix estimation under {H}uber's
  contamination model.
\newblock {\em Ann. Statist.}, 46(5):1932--1960, 2018.

\bibitem{MR3289373}
Yanchu Chen and Qi-Man Shao.
\newblock Berry-{E}sseen inequality for unbounded exchangeable pairs.
\newblock In {\em Probability approximations and beyond}, volume 205 of {\em
  Lect. Notes Stat.}, pages 13--30. Springer, New York, 2012.

\bibitem{Bartlett19}
Yeshwanth Cherapanamjeri, Nicolas Flammarion, and Peter~L. Bartlett.
\newblock Fast mean estimation with sub-gaussian rates, 2019.

\bibitem{dalalyan2019outlier}
Arnak Dalalyan and Philip Thompson.
\newblock Outlier-robust estimation of a sparse linear model using l1-penalized
  huber's m-estimator.
\newblock In {\em Advances in Neural Information Processing Systems}, pages
  13188--13198, 2019.

\bibitem{dalalyan2020all}
Arnak~S Dalalyan and Arshak Minasyan.
\newblock All-in-one robust estimator of the gaussian mean.
\newblock {\em arXiv preprint arXiv:2002.01432}, 2020.

\bibitem{MR902258}
P.~L. Davies.
\newblock Asymptotic behaviour of {$S$}-estimates of multivariate location
  parameters and dispersion matrices.
\newblock {\em Ann. Statist.}, 15(3):1269--1292, 1987.

\bibitem{MR2752147}
Michiel Debruyne.
\newblock An outlier map for support vector machine classification.
\newblock {\em Ann. Appl. Stat.}, 3(4):1566--1580, 2009.

\bibitem{Jules_Guillaume_1}
Jules Depersin and Guillaume Lecu{\'e}.
\newblock Fast algorithms for robust estimation of a mean vector.
\newblock 2019.

\bibitem{convex_prog}
Jules Depersin and Guillaume Lecu{\'e}.
\newblock Convex programs and algorithms for robust subgaussian estimation of a
  mean vector with respect to any norm.
\newblock Technical report, IPParis, Crest, ENSAE, 2020.

\bibitem{devroye2016}
Luc Devroye, Matthieu Lerasle, Gabor Lugosi, and Roberto~I. Oliveira.
\newblock Sub-gaussian mean estimators.
\newblock {\em Ann. Statist.}, 44(6):2695--2725, 12 2016.

\bibitem{MR3576558}
Luc Devroye, Matthieu Lerasle, Gabor Lugosi, and Roberto~I. Oliveira.
\newblock Sub-{G}aussian mean estimators.
\newblock {\em Ann. Statist.}, 44(6):2695--2725, 2016.

\bibitem{MR3631028}
Ilias Diakonikolas, Gautam Kamath, Daniel~M. Kane, Jerry Li, Ankur Moitra, and
  Alistair Stewart.
\newblock Robust estimators in high dimensions without the computational
  intractability.
\newblock In {\em 57th {A}nnual {IEEE} {S}ymposium on {F}oundations of
  {C}omputer {S}cience---{FOCS} 2016}, pages 655--664. IEEE Computer Soc., Los
  Alamitos, CA, 2016.

\bibitem{diakonikolas2017being}
Ilias Diakonikolas, Gautam Kamath, Daniel~M Kane, Jerry Li, Ankur Moitra, and
  Alistair Stewart.
\newblock Being robust (in high dimensions) can be practical.
\newblock {\em arXiv preprint arXiv:1703.00893}, 2017.

\bibitem{donoho1982}
D~Donoho.
\newblock {\em L.(1982) Breakdown properties of multivariate location
  estimtors}.
\newblock PhD thesis, Ph. D. Qualifying Paper, Dept. of Statistics, Harvard
  Univ.

\bibitem{MR689745}
David Donoho and Peter~J. Huber.
\newblock The notion of breakdown point.
\newblock In {\em A {F}estschrift for {E}rich {L}. {L}ehmann}, Wadsworth
  Statist./Probab. Ser., pages 157--184. Wadsworth, Belmont, CA, 1983.

\bibitem{donoho1982breakdown}
David~L Donoho.
\newblock Breakdown properties of multivariate location estimators.
\newblock Technical report, Technical report, Harvard University, Boston. URL
  http://www-stat. stanford~…, 1982.

\bibitem{MR1193313}
David~L. Donoho and Miriam Gasko.
\newblock Breakdown properties of location estimates based on halfspace depth
  and projected outlyingness.
\newblock {\em Ann. Statist.}, 20(4):1803--1827, 1992.

\bibitem{donoho1992breakdown}
David~L Donoho and Miriam Gasko.
\newblock Breakdown properties of location estimates based on halfspace depth
  and projected outlyingness.
\newblock {\em The Annals of Statistics}, 20(4):1803--1827, 1992.

\bibitem{gnanadesikan1972robust}
Ramanathan Gnanadesikan and John~R Kettenring.
\newblock Robust estimates, residuals, and outlier detection with multiresponse
  data.
\newblock {\em Biometrics}, pages 81--124, 1972.

\bibitem{haldane1948note}
JBS Haldane.
\newblock Note on the median of a multivariate distribution.
\newblock {\em Biometrika}, 35(3-4):414--417, 1948.

\bibitem{MR0359096}
Frank~R. Hampel.
\newblock Robust estimation: a condensed partial survey.
\newblock {\em Z. Wahrscheinlichkeitstheorie und Verw. Gebiete}, 27:87--104,
  1973.

\bibitem{MR362657}
Frank~R. Hampel.
\newblock The influence curve and its role in robust estimation.
\newblock {\em J. Amer. Statist. Assoc.}, 69:383--393, 1974.

\bibitem{hopkins2018sub}
Samuel~B Hopkins.
\newblock Sub-gaussian mean estimation in polynomial time.
\newblock {\em arXiv preprint arXiv:1809.07425}, 2018.

\bibitem{MR2488795}
Peter~J. Huber and Elvezio~M. Ronchetti.
\newblock {\em Robust statistics}.
\newblock Wiley Series in Probability and Statistics. John Wiley \& Sons, Inc.,
  Hoboken, NJ, second edition, 2009.

\bibitem{MR855970}
Mark~R. Jerrum, Leslie~G. Valiant, and Vijay~V. Vazirani.
\newblock Random generation of combinatorial structures from a uniform
  distribution.
\newblock {\em Theoret. Comput. Sci.}, 43(2-3):169--188, 1986.

\bibitem{MR3136463}
B.~Klartag and S.~Sodin.
\newblock Variations on the {B}erry-{E}sseen theorem.
\newblock {\em Teor. Veroyatn. Primen.}, 56(3):514--533, 2011.

\bibitem{MR2520120}
Bo'az Klartag.
\newblock A {B}erry-{E}sseen type inequality for convex bodies with an
  unconditional basis.
\newblock {\em Probab. Theory Related Fields}, 145(1-2):1--33, 2009.

\bibitem{MR3645123}
Bo'az Klartag.
\newblock Super-{G}aussian directions of random vectors.
\newblock In {\em Geometric aspects of functional analysis}, volume 2169 of
  {\em Lecture Notes in Math.}, pages 187--211. Springer, Cham, 2017.

\bibitem{Kol11}
V.~Koltchinskii.
\newblock {\em Oracle {I}nequalities in {E}mpirical {R}isk {M}inimization and
  {S}parse {R}ecovery {P}roblems}.
\newblock Springer, Berlin, 2011.

\bibitem{MR4102681}
Guillaume Lecu\'{e} and Matthieu Lerasle.
\newblock Robust machine learning by median-of-means: theory and practice.
\newblock {\em Ann. Statist.}, 48(2):906--931, 2020.

\bibitem{MR2814399}
Michel Ledoux and Michel Talagrand.
\newblock {\em Probability in {B}anach spaces}.
\newblock Classics in Mathematics. Springer-Verlag, Berlin, 2011.
\newblock Isoperimetry and processes, Reprint of the 1991 edition.

\bibitem{lei2020fast}
Zhixian Lei, Kyle Luh, Prayaag Venkat, and Fred Zhang.
\newblock A fast spectral algorithm for mean estimation with sub-gaussian
  rates.
\newblock In {\em Conference on Learning Theory}, pages 2598--2612. PMLR, 2020.

\bibitem{MR1091202}
O.~V. Lepski\u{\i}.
\newblock A problem of adaptive estimation in {G}aussian white noise.
\newblock {\em Teor. Veroyatnost. i Primenen.}, 35(3):459--470, 1990.

\bibitem{MR1147167}
O.~V. Lepski\u{\i}.
\newblock Asymptotically minimax adaptive estimation. {I}. {U}pper bounds.
  {O}ptimally adaptive estimates.
\newblock {\em Teor. Veroyatnost. i Primenen.}, 36(4):645--659, 1991.

\bibitem{LO}
M.~Lerasle and R.~Oliveira.
\newblock Robust empirical mean estimators.
\newblock Technical report, IMPA and CNRS, 2011.

\bibitem{lerasle2019monk}
Matthieu Lerasle, Zolt{\'a}n Szab{\'o}, Timoth{\'e}e Mathieu, and Guillaume
  Lecu{\'e}.
\newblock Monk outlier-robust mean embedding estimation by median-of-means.
\newblock In {\em International Conference on Machine Learning}, pages
  3782--3793. PMLR, 2019.

\bibitem{MR1041400}
Regina~Y. Liu.
\newblock On a notion of data depth based on random simplices.
\newblock {\em Ann. Statist.}, 18(1):405--414, 1990.

\bibitem{MR1214839}
Regina~Y. Liu.
\newblock Data depth and multivariate rank tests.
\newblock In {\em {$L_1$}-statistical analysis and related methods
  ({N}euch\^{a}tel, 1992)}, pages 279--294. North-Holland, Amsterdam, 1992.

\bibitem{MR1186260}
Regina~Y. Liu and Kesar Singh.
\newblock Ordering directional data: concepts of data depth on circles and
  spheres.
\newblock {\em Ann. Statist.}, 20(3):1468--1484, 1992.

\bibitem{MR3217437}
Karim Lounici.
\newblock High-dimensional covariance matrix estimation with missing
  observations.
\newblock {\em Bernoulli}, 20(3):1029--1058, 2014.

\bibitem{lu2020robust}
Junwei Lu, Fang Han, and Han Liu.
\newblock Robust scatter matrix estimation for high dimensional distributions
  with heavy tail.
\newblock {\em IEEE transactions on information theory}, 2020.

\bibitem{MR4017683}
G\'{a}bor Lugosi and Shahar Mendelson.
\newblock Mean estimation and regression under heavy-tailed distributions: a
  survey.
\newblock {\em Found. Comput. Math.}, 19(5):1145--1190, 2019.

\bibitem{MR4026610}
G\'{a}bor Lugosi and Shahar Mendelson.
\newblock Near-optimal mean estimators with respect to general norms.
\newblock {\em Probab. Theory Related Fields}, 175(3-4):957--973, 2019.

\bibitem{lugosi2019sub}
G{\'a}bor Lugosi, Shahar Mendelson, et~al.
\newblock Sub-gaussian estimators of the mean of a random vector.
\newblock {\em The Annals of Statistics}, 47(2):783--794, 2019.

\bibitem{LMSL}
Z.~Szabo M.~Lerasle, T.~Matthieu and G.~Lecué.
\newblock Monk – outliers-robust mean embedding estimation by
  median-of-means.
\newblock Technical report, CNRS, University of Paris 11, Ecole Polytechnique
  and CREST, 2017.

\bibitem{MR2238141}
Ricardo~A. Maronna, R.~Douglas Martin, and Victor~J. Yohai.
\newblock {\em Robust statistics}.
\newblock Wiley Series in Probability and Statistics. John Wiley \& Sons, Ltd.,
  Chichester, 2006.
\newblock Theory and methods.

\bibitem{maronna1995behavior}
Ricardo~A Maronna and Victor~J Yohai.
\newblock The behavior of the stahel-donoho robust multivariate estimator.
\newblock {\em Journal of the American Statistical Association},
  90(429):330--341, 1995.

\bibitem{MS}
S~Minsker and N.~Strawn.
\newblock Distributed statistical estimation and rates of convergence in normal
  approximation.
\newblock Technical report, arXiv: 1704.02658, 2017.

\bibitem{minsker2015geometric}
Stanislav Minsker.
\newblock Geometric median and robust estimation in banach spaces.
\newblock {\em Bernoulli}, 21(4):2308--2335, 2015.

\bibitem{minsker2018uniform}
Stanislav Minsker.
\newblock Uniform bounds for robust mean estimators.
\newblock {\em arXiv preprint arXiv:1812.03523}, 2018.

\bibitem{nagy2019halfspace}
Stanislav Nagy, Carsten Sch{\"u}tt, Elisabeth~M Werner, et~al.
\newblock Halfspace depth and floating body.
\newblock {\em Statistics Surveys}, 13:52--118, 2019.

\bibitem{MR702836}
A.~S. Nemirovsky and D.~B.~and Yudin.
\newblock {\em Problem complexity and method efficiency in optimization}.
\newblock A Wiley-Interscience Publication. John Wiley \& Sons, Inc., New York,
  1983.
\newblock Translated from the Russian and with a preface by E. R. Dawson,
  Wiley-Interscience Series in Discrete Mathematics.

\bibitem{pena2007combining}
Daniel Pe{\~n}a and Francisco~J Prieto.
\newblock Combining random and specific directions for outlier detection and
  robust estimation in high-dimensional multivariate data.
\newblock {\em Journal of Computational and Graphical Statistics},
  16(1):228--254, 2007.

\bibitem{MR1353441}
Valentin~V. Petrov.
\newblock {\em Limit theorems of probability theory}, volume~4 of {\em Oxford
  Studies in Probability}.
\newblock The Clarendon Press, Oxford University Press, New York, 1995.
\newblock Sequences of independent random variables, Oxford Science
  Publications.

\bibitem{rousseeuw2018measure}
Peter~J Rousseeuw, Jakob Raymaekers, and Mia Hubert.
\newblock A measure of directional outlyingness with applications to image data
  and video.
\newblock {\em Journal of Computational and Graphical Statistics},
  27(2):345--359, 2018.

\bibitem{stahel1981robuste}
Werner~A Stahel.
\newblock {\em Robuste sch{\"a}tzungen: infinitesimale optimalit{\"a}t und
  sch{\"a}tzungen von kovarianzmatrizen}.
\newblock PhD thesis, ETH Zurich, 1981.

\bibitem{MR0426989}
John~W. Tukey.
\newblock Mathematics and the picturing of data.
\newblock In {\em Proceedings of the {I}nternational {C}ongress of
  {M}athematicians ({V}ancouver, {B}. {C}., 1974), {V}ol. 2}, pages 523--531,
  1975.

\bibitem{MR1292555}
David~E. Tyler.
\newblock Finite sample breakdown points of projection based multivariate
  location and scatter statistics.
\newblock {\em Ann. Statist.}, 22(2):1024--1044, 1994.

\bibitem{van2016stahel}
Stefan Van~Aelst.
\newblock Stahel--donoho estimation for high-dimensional data.
\newblock {\em International Journal of Computer Mathematics}, 93(4):628--639,
  2016.

\bibitem{van2011stahel}
Stefan Van~Aelst, E~Vandervieren, and Gert Willems.
\newblock Stahel-donoho estimators with cellwise weights.
\newblock {\em Journal of Statistical Computation and Simulation}, 81(1):1--27,
  2011.

\bibitem{MR1385671}
Aad~W. van~der Vaart and Jon~A. Wellner.
\newblock {\em Weak convergence and empirical processes}.
\newblock Springer Series in Statistics. Springer-Verlag, New York, 1996.
\newblock With applications to statistics.

\bibitem{MR3408730}
V.~N. Vapnik and A.~Ya. Chervonenkis.
\newblock On the uniform convergence of relative frequencies of events to their
  probabilities.
\newblock In {\em Measures of complexity}, pages 11--30. Springer, Cham, 2015.
\newblock Reprint of Theor. Probability Appl. {{\bf{1}}6} (1971), 264--280.

\bibitem{MR1719582}
Vladimir~N. Vapnik.
\newblock {\em The nature of statistical learning theory}.
\newblock Statistics for Engineering and Information Science. Springer-Verlag,
  New York, second edition, 2000.

\bibitem{MR2051003}
Yijun Zuo, Hengjian Cui, and Xuming He.
\newblock On the {S}tahel-{D}onoho estimator and depth-weighted means of
  multivariate data.
\newblock {\em Ann. Statist.}, 32(1):167--188, 2004.

\end{thebibliography}
\end{footnotesize}

\end{document}